\documentclass{article}

\usepackage{color}

\newcommand{\ssh}[1]{{\textcolor{red}{#1}}}

\usepackage{amsfonts}
\usepackage{amsmath,amsthm,amssymb}

\usepackage{latexsym}
\usepackage{enumerate}

\newcommand\bbZ{{\mathbb Z}}

\newcommand{\bbR}{{\mathbb{R}}}

\newcommand\oZ{\operatorname{Z}}
\newcommand\oC{\operatorname{C}}
\newcommand\oP{\operatorname{P}}

\newcommand\oO{\operatorname{O}}

\newcommand\lla{{\langle}}
\newcommand\rra{{\rangle}}

\newtheorem{thm}{Theorem}[section]
\newtheorem{prop}[thm]{Proposition}

\newtheorem{lem}[thm]{Lemma}
\newtheorem{cor}[thm]{Corollary}

\numberwithin{equation}{section}

\newcommand\pf{\noindent{\bf Proof.}~~}

\allowdisplaybreaks[1]

\title{The spectra of finite $3$-transposition groups}
\author{J.I. Hall\\
	Department of Mathematics\\
	Michigan State University\\
	619 Red Cedar Road\\ 
	East Lansing MI 48824 USA\\
	\tt{halljo@msu.edu}
	\and 
	S.~Shpectorov\\
	School of Mathematics\\
	Watson Building\\
	University of Birmingham\\
	Edgbaston\\
	Birmingham B15 2TT UK\\
	\tt{s.shpectorov@bham.ac.uk}
}

\begin{document}
\maketitle

{\abstract{We calculate the spectrum of the diagram
		for each finite $3$-transposition group. 
		Such graphs with a given minimum eigenvalue have
		occurred in the context of compact Griess subalgebras
		of vertex operator algebras.}}

\hyphenation{trans-pos-i-tion}
\hyphenation{trans-pos-i-tions}

\newcommand{\dspace}[1]{\vspace*{#1\baselineskip}}
\newcommand{\trs}{trans\-pos\-i\-tion}
\newcommand{\trans}{$3$-trans\-pos\-i\-tion}
\newcommand{\su}[2][2]{\operatorname{SU}_{#2}(#1)}
\newcommand{\spt}[2][2]{\operatorname{Sp}_{#2}(#1)}
\newcommand{\tri}[1]{\operatorname{P}\!\Omega_{8}^{+}(#1)\splt\Sym{3}}
\newcommand{\eps}{\epsilon}
\newcommand{\lam}{\lambda}
\newcommand{\ogpt}[2][\eps]{\operatorname{O}_{#2}^{#1}(2)}
\newcommand{\gp}[6]{\mbox{}^{#4}_{#5}{#1}_{#2}^{#3}{#6}}
\newcommand{\omgp}[2]{\gp{\Omega}{#1}{#2}{+}{}{(3)}}
\newcommand{\ogp}[2]{\gp{\operatorname{O}}{#1}{#2}{+}{}{(3)}}
\newcommand{\igp}[1]{\gp{\Omega}{}{}{+}{}{#1}}
\newcommand{\ofive}{\omgp{5}{-}}
\newcommand{\tosix}{\mbox{${3}\nsplt\omgp{6}{-}$}}
\newcommand{\toseven}{\mbox{${3}\nsplt\omgp{7}{+}$}}
\newcommand{\ossix}{{}^{+}\Omega_{6}^{-}(3)}
\newcommand{\osix}{\mbox{${3} \nsplt\ossix$}}
\newcommand{\Fi}[1]{\operatorname{Fi}_{{#1}}}
\newcommand{\codiag}[1]{[#1]}
\newcommand{\diag}[1]{(#1)}
\newcommand{\name}[2][0]{\dspace{#1}{\rm\bfseries{#2}}}
\newcommand{\comd}[1]{{\bullet #1}}
\newcommand{\splt}{\colon\!\!}
\newcommand{\nsplt}{\,{}^{\displaystyle .}}
\newcommand{\allone}{\mathbf{1}}
\newcommand{\Spec}[1]{\operatorname{Spec}(#1)}
\newcommand{\spec}[1]{(\!(#1)\!)}
\newcommand{\AMat}[1]{\operatorname{AMat}\!\left(#1\right)}
\newcommand{\Sym}[1]{\operatorname{Sym}(#1)}
\newcommand{\Weyl}[2][]{\operatorname{W}_{#1}(#2)}
\newcommand{\iso}{\simeq}
\newcommand{\gf}[1]{{\mathbb F}_{#1}}
\newcommand{\announce}[1]{\emph{#1}}
\newcommand{\adj}{\,\sim\,}
\newcommand{\adjns}{\sim}
\newcommand{\nadj}{\,\not\sim\,}
\newcommand{\eigmult}[2]{[#1]^{#2}}
\renewcommand{\omgp}[3][]{\gp{\Omega}{#2}{#3}{+}{\,#1}{(3)}}
\newcommand{\presub}[1]{\mbox{}_{#1}}
\newcommand{\onehalf}{\frac{1}{2}}

\section{Introduction}

A $3$-transposition group $(G,D)$ is a group $G$ generated
by a conjugacy class $D=D^G$ of elements of order $2$
such that
\[
d,e \in D \implies |de| \in \{1,2,3\}\,.
\]
The \announce{diagram} of $(G,D)$ is the graph with vertex set
$D$ and edges
\[
d \adj e \iff |de|=3\,.
\]
As a consequence of work by Fischer \cite{Fi71} and 
later Cuypers and Hall
\cite{CuHa95} all diagrams for all finite $3$-transposition
groups are known. In this paper we give the eigenvalues
and spectrum of (the adjacency matrix of) each such diagram.

Of particular importance are the minimum eigenvalues, always
a negative integer. One result is that for given 
$-t$ the possibilities for $3$-transposition groups 
with minimum eigenvalue (greater than or) equal to $-t$
are limited. 
Miyamoto \cite{Mi96} first observed a connection
between $3$-transposition groups and
compact Griess subalgebras found within vertex
operator algebras. Particularly relevant for this paper is 
the work of Matsuo \cite{Ma03,Ma05}.

The complement of the diagram is the \announce{codiagram}
or \announce{commuting graph}.

\medskip
{\bf Acknowledgement.} The authors are very grateful to the 
referees of this paper for careful reading, correcting a few 
misprints and errors, and generally helping to make this rather 
technical paper more readable. 

\section{Eigenvalues of graphs}

Let $X$ be a nonempty set and $\diag{X}$ a graph with $X$ as vertex set.
The adjacency of the edge $(x,y)$ is written $x \adj y$. The $(0,1)$-adjacency matrix
of the graph will be denoted $\AMat{\diag{X}}$, and the spectrum of the graph
is the (ordered) spectrum of $\AMat{\diag{X}}$:
\[
  \Spec{\diag{X}}= \spec{\dots,r_i, \dots}\,.
\]

The all-one vector $\allone$ is an eigenvector of $\AMat{\diag{X}}$ with eigenvalue 
$k$ if and only if $\diag{X}$ is regular of degree $k$. 
In this case, by the Perron-Frobenius Theorem, $k$ is 
    the largest eigenvalue and the corresponding eigenspace has dimension the number of 
    connected components of $\diag{X}$. Hence, when $\diag{X}$ is connected, which will 
    (almost) always be the case in this paper, the eigenspace for $k$ is $1$-dimensional, 
    spanned by $\allone$.

Furthermore, 
all other eigenspaces of the regular 
    connected graph $\diag{X}$ are perpendicular to $\allone$; that is, they belong to 
    the sum-zero hyperplane of $\bbR^n$. Such eigenvectors and their associated 
    eigenvalues will be called \announce{restricted}.

For this reason, we will list $k$ first in the spectrum and separate it from the 
    restricted eigenvalues by a semicolon.

The complement of the graph $\diag{X}$ is the graph $\codiag{X}$
with the same vertex set but all edges replaced by nonedges and
nonedges by edges. Thus
\[
  \AMat{\codiag{X}}=J_n -I_n - \AMat{\diag{X}}\,.
\]
where $n=|X|$, $I_n$ is the $n \times n$ identity matrix,
and $J_n$ is the $n \times n$ matrix consisting entirely of 
$1$'s. (We may drop the subscripts, when $n$ is apparent.)
All nonzero vectors of $\bbR^n$ are eigenvectors of $I_n$ with eigenvalue
$1$. The all-one vector $\allone$ is an eigenvector of $J_n$ with eigenvalue
$n$ of multiplicity $1$, and the sum-zero hyperplane of $\bbR^n$ consists 
of null vectors for $J_n$---its nonzero vectors are eigenvectors with 
eigenvalue $0$. 

We thus have

\begin{prop}\label{prop-comp-eig}
If $\diag{X}$ is a regular graph of degree $k$ and 
the spectrum of $\diag{X}$ is $\spec{k; \dots,r_i,\dots}$,
then the spectrum of $\codiag{X}$ is $\spec{l; \dots,-1-r_i,\dots}$
where $|X|=1+k+l$.
\qed
\end{prop}

If $M$ is an $n \times n$ matrix, then $2^\comd{1}M$ is the $2n \times 2n$ matrix
\[
 \begin{bmatrix}
  M & M \\ M & M
 \end{bmatrix}
 = M \otimes J_2
\,;
\]
and $2^\comd{h}M$ is the $2^hn \times 2^hn$ matrix that results from 
repeating this construction $h$ times.

If $M$ is an $n \times n$ matrix, then $3^\comd{1}M$ is the $3n \times 3n$ matrix
\[
 \begin{bmatrix}
  M & M +I_n & M+I_n \\ M+I_n & M & M+I_n\\M+I_n&M+I_n&M
 \end{bmatrix}
 = (M+I_n) \otimes J_3 -I_{3n}
\,;
\]
and $3^\comd{h}M$ is the $3^hn \times 3^hn$ matrix that results from 
repeating this construction $h$ times.

This nonstandard notation $p^{\comd{h}}M$ and the relevance of these
	matrices will become clear with Lemma \ref{lem-adj-mat} and the remarks
	that proceed it.

\begin{prop}
\label{prop-eig23}
Let $v_1(=\allone), \dots, v_i,\dots, v_n$ be a basis of eigenvectors for the matrix
$M$, the associated spectrum being $\spec{\dots,r_i,\dots}$.
\begin{enumerate}[{\rm(a)}]
 \item $2^{\comd{1}}M$ has the basis of eigenvectors
\[
(v_1,v_1)\,,\, \dots\,,\, (v_i,v_i)\,,\, (v_i,-v_i)\,,\, \dots\,,\, (v_n,-v_n)
\]
with associated
spectrum $\spec{\dots, 2r_i,0, \dots}$.
 \item
$3^{\comd{1}}M$ has the basis of eigenvectors
\begin{multline*}
 (v_1,v_1,v_1)\,,\, (v_1,-v_1,0)\,,\,\dots\,, (v_i,v_i,v_i)\,,\, (v_i,-v_i,0)\,,\,
 (v_i,0,-v_i)\,,\\ \dots\,,\,(v_n,-v_n,0)\,,\,(v_n,0,-v_n)
\end{multline*}
with associated spectrum $\spec{\dots, 3r_i+2,-1,-1 \dots}$. \qed
\end{enumerate}
\end{prop}

If $M$ is an $n \times n$ matrix, then $3 \times M$ is the $3n \times 3n$ matrix
\[
 \begin{bmatrix}
  M & J_n & J_n \\ J_n & M & J_n\\J_n&J_n&M
 \end{bmatrix}
\,.
\]
As the restricted eigenvectors of a regular graph are also eigenvectors
for $J_n$, we find

\begin{prop}\label{prop-eigtri}
Let $v_1(=\allone), \dots, v_i,\dots, v_n$ be a basis of eigenvectors for the
adjacency matrix $M$ of a regular graph of degree $k$,
the associated spectrum being $\spec{k; \dots,r_i,\dots}$.
Then $3\times M$ has the basis of eigenvectors
\begin{multline*}
 (v_1,v_1,v_1), (v_1,-v_1,0), (v_1,0,-v_1),
\,\dots\,, (v_i,v_i,v_i), (v_i,-v_i,0),  (v_i,0,-v_i),\dots\\
\dots\,, (v_n,v_n,v_n), (v_n,-v_n,0),  (v_n,0,-v_n)
\end{multline*}
with associated spectrum 
$\spec{k+2n; -(n-k), -(n-k),\, \dots\,,\, r_i,r_i,r_i,\, \dots}$.
In particular, $-(n-k)$ is an eigenvalue of $3 \times M$. \qed
\end{prop}

The nonstandard matrix notation $3 \times M$ and these
	results will reappear in Theorem \ref{thm-tri}.

\section{Rank $\mathbf 3$ and strongly regular graphs}
\label{sec-rk3}

Consider a graph $\diag{X}$ and subgroup $G$ of its automorphism group with the 
following property:
\begin{quote}
$G$ is transitive on $X$, on the set of ordered edges of $\diag{X}$, and
on the set of ordered edges of $\codiag{X}$.
\end{quote}
Assuming that all three sets are nonempty, we say that $G$
acts with \announce{rank $3$} on $\diag{X}$ (and so also on $\codiag{X}$)
and that $\diag{X}$ and $\codiag{X}$ are a complementary pair
of \announce{rank $3$ graphs}. 
(There is nothing to say if all three are empty. If two are empty,
then $|X|=1$. If one is empty, 
then $\diag{X}$ and $\codiag{X}$ are a complementary pair of a complete
and an empty graph, and $G$ is $2$-transitive on $X$; this is 
\announce{rank $2$} action.)

A \announce{strongly regular graph} is a finite graph $\diag{X}$ with
the following strong regularity property:
\begin{quote}
There are constants $k$, $\lam$, and $\mu$ such that,
for $x,y \in X$, the number of common neighbors of $x,y$ is
$k$ when $x = y$; $\lambda$ when $x \adj y$; and $\mu$ when $x \nadj y$.
\end{quote}
Empty and complete graphs provide the
degenerate cases $k=0$ and $k=n-1$ of this condition,
where
\begin{equation}
	|X|=n\,.
\end{equation}
Here we do not include these as strongly regular; that is,
we additionally require 
\begin{equation}
	0 < k < n-1\,.
\end{equation}
This graph will be connected of diameter $2$ unless $\mu=0$. In that case, the graph 
is a disjoint union	of complete subgraphs $K_{k+1}$. Its complementary strongly regular 
graph is then complete multipartite with $\mu=k$. This pair of graphs is 
\announce{imprimitive}. We shall only be concerned with strongly regular graphs that are 
not imprimitive---those that are \announce{primitive}.

For us, the basic 
observation is that a rank $3$ graph is strongly regular.
A strongly regular graph is, in particular, regular of degree $k$.
One says that the strongly regular graph $\diag{X}$ has \announce{parameters}
$(n,k,\lambda,\mu)$.
The parameters are thus nonnegative integers with 
\begin{equation}
	n > k \ge \mu
	\quad\text{and}\quad 
	k-1\ge \lambda\,.
\end{equation}

An elementary calculation shows that if $\diag{X}$ is strongly
regular with parameters $(n,k,\lambda,\mu)$, then $\codiag{X}$
is also strongly regular, its parameters being
\begin{equation}\label{eqn-pp}
(n,k',\lambda',\mu') 
\quad\text{for}\quad 
k'=n-k-1\,,\
\lam'=n-2k+\mu-2\,, 
\ \text{and}\ 
\mu'=n-2k+\lambda\,.
\end{equation}

It is usual to write the codegree as
\begin{equation}
	l=k'=n-k-1\,.
\end{equation}
By counting paths of length two from a fixed
vertex
\begin{equation} \label{eqn-mul}
	\mu l = k(k-1-\lambda) \,.
\end{equation}

\dspace{.5}

Let $M$ be the adjacency matrix of the strongly regular graph
$\diag{X}$ with parameters $(n,k,\lambda,\mu)$. Counting
all directed paths of length $2$ yields
\[
  M^2= kI + \lam M + \mu(J-I-M)
\]
whence
\[
  M^2 + (\mu-\lam)M + (\mu-k)I = \mu J\,.
\]
In particular, the restricted eigenvalues of $M$ are the roots
$r$ and $s$ of the monic quadratic polynomial 
\begin{equation}
	  x^2 + (\mu-\lam)x + (\mu-k)\,.
\end{equation}
As $\mu \le k$ the roots $r$ and $s$ are real.
We take $s \le r$ by convention. As ${\rm tr}\,M = 0$ and $k$ is a 
positive eigenvalue of $A$, $s<0$. Also, $-rs=k-\mu\geq 0$, and so 
the real parameters $s$ and $r$ are restricted by
\begin{equation}
	s<0 \text{ and }0 \le r \le k\,.
\end{equation}
In particular, $s<r$.

The triple $(k,r,s)$ is determined by $(k,\lambda,\mu)$.
Conversely $(k,r,s)$ determines $(k,\lambda,\mu)$ via
\begin{equation}
	\mu = k+rs \quad\text{and}\quad \lambda = \mu + r + s\,.
\end{equation}

Let $r$ and $s$  have (restricted) multiplicities $f$ and $g$,
respectively. As $r \neq s$, the parameters
$f,g$ can be found from 
\begin{equation} 
	1+f+g = n 
	\quad\text{and}\quad 
	k+fr+gs = {\rm tr}\,M = 0\,.
\end{equation}
Indeed the multiplicities are
\begin{equation}
	f=(r-s)^{-1}(-sn+s-k)
	\quad\text{and}\quad
	g=(r-s)^{-1}(rn-r+k)\,.
\end{equation}
The fact that $f,g$ must be integers is a strong restriction
on possible parameter sets.

Conversely, given the integer $f$ and $g$,
if $f = g$, then $f = g = (n-1)/2$. Therefore
\[
  k =-fr-gs=-(r+s)(n-1)/2=(\mu-\lambda)(n-1)/2\,.
\]
Since $0 < k < n-1$,
it follows that 
$\mu = \lambda+1\neq 0$ and $k = (n-1)/2=l$.
As $k=l$, we have $\mu = k-1-\lambda$, hence $k=2\mu$. Thus
$(n,k,\lambda,\mu) = (4t+1,2t,t-1,t)$ for a suitable integer $t$,
and $r,s = (-1\pm\sqrt{n})/2$. 
This is known as the \announce{half case}
and will not be of concern here.

In the generic case $f \ne g$, one can solve for $r,s$ from 
\begin{equation}
	r+s=\lambda-\mu
	\quad\text{and}\quad 
	fr+gs = -k\,. 
\end{equation}
It follows that $r,s$ are rational.
As roots of a monic polynomial with integer coefficients,
they are also algebraic integers; so they are integral in this case.

The \announce{extended parameter} list for $\diag{X}$ is
\[
  (n,k,\lam,\mu\,;\eigmult{r}{f},\eigmult{s}{g})
\]
or
\[
  (n,k,\lam,\mu\,;\{\eigmult{r}{f},\eigmult{s}{g}\})
\]
when it is not
clear which eigenvalue is $r$ and which is $s$. The corresponding
extended parameter list for $\codiag{X}$ is
\[
(n,l,\lam',\mu'\,;\eigmult{r'}{g},\eigmult{s'}{f})
\]
or
\[
(n,l,\lam',\mu'\,;\{\eigmult{r'}{g},\eigmult{s'}{f}\})
\]
where
\begin{equation}
	\label{eqn-prime-param}
	l=n-k-1\,,\
	\lam'=n-2k+\mu-2\,, 
	\mu'=n-2k+\lambda
\end{equation}
as before, and additionally
\begin{equation}
r'=-s-1\,,\quad 
s'=-r-1\,
\end{equation}
since $\AMat{\codiag{X}}+\AMat{\diag{X}}=J_n -I_n $;
see Proposition \ref{prop-comp-eig}.

\dspace{.5}
These sets of parameters
are highly redundant, being related by the various equations
of this section. All parameters can be determined
by various small subsets of the complete parameter list. 
In particular three parameters are enough when we have
\[
n\,;\quad \text{one of}\  \
k=l'\ \text{or}\ l = k'\,;
\quad\text{any one of}\ \ \lam,\mu,\lam',\mu'\,.
\]
Of course, the more parameters that can be calculated
directly, the easier the remaining calculations will be.

It is also of note that all parameters can be
derived from the spectrum
\[
\spec{k;\eigmult{r}{f},\eigmult{s}{g}}\,.
\]

We have already seen that the values $\mu=0$ and $\mu=k$
are special---these are the imprimitive graphs.
Indeed these parameters make the complementary statements 
that one of $\diag{X}$ or $\codiag{X}$ is a nontrivial equivalence 
relation---a disjoint union of complete subgraphs (of fixed size 
$m>1$)---while the other is  a complete multipartite graph with all 
parts of size $m$.
As an important special case, when $G$ acts imprimitively with rank $3$ 
on $\diag{X}$ and $\codiag{X}$, these form a complementary pair of 
imprimitive strongly regular graphs.

\section{{$\mathbf 3$}-transposition diagrams and eigenvalues}

The normal set $D$ of the group $G$ is a set of 
\announce{$3$-transpositions} in $G$ if it consists of
elements of order $2$ with the property
\[
d,e \in D \implies |de| \in \{1,2,3\}\,.
\]
The study of such sets $D$ and groups $G$ was initiated
by Bernd Fischer \cite{Fi71}. Fischer's paper and
the later paper \cite{CuHa95} of Cuypers and Hall are
our basic references on this topic.

If $E$ is a subset of $D$ in $G$ then the
\announce{diagram} of $E$, denoted $\diag{E}$, is the graph with vertex
set $E$ and having an edge between the two vertices $d,e$
precisely when $|de|=3$. The \announce{commuting graph} of $E$,
or \announce{codiagram} of $E$, is the graph complement $\codiag{E}$ 
of the diagram of $E$.

There are two cases of primary interest. The first has
$E$ some small generating set of $G$; for
instance, the \trans\ group $\Sym{n+1}$ is the Weyl group $\Weyl{A_n}$
with diagram the $n$-vertex path $A_n$. In the second case $E$ is 
equal to the full
class $D$, and we then abuse terminology by
saying that the diagram of $D$ is also the diagram of $G$.

\begin{thm}\label{thm-trans-prop}
	\mbox{}

	\begin{enumerate}[$(a)$]
		\item If $H$ is a subgroup of $G$, then 
		$D \cap H=\emptyset$ or
		$D \cap H$ is a normal 
        set of $3$-transpositions in $H$. If $N$ is a normal subgroup 
        of $G$, then 
        $D \subset N$ or the nontrivial elements of $DN/N$ form
        a normal subset of $3$-transpositions in $G/N$.
        \item
        Let $D_i$, for $i \in I$, be the connected components of
        $\diag{D}$. Then each $D_i$ is a conjugacy class of 
	    $3$-transpositions in the group $G_i=\langle D_i \rangle$.
	    Furthermore the normal subgroup $\langle D \rangle$ of $G$
	    is the central product of its subgroups $G_i$.
	    \item
	    If $G=\langle D \rangle$ then, for $d \in D \backslash \oZ(G)$, 
	    each coset $d \oZ (G)$ meets $D$ only in $d$. 
	\end{enumerate}
\end{thm}

The first two parts of the theorem are Fischer's basic
Inheritance Properties \cite[(1.2)]{Fi71}. The second of
these allows us to focus on the case
$G = \langle D \rangle$ for the conjugacy class $D$ of
$3$-transpositions.
In this situation we say that $(G,D)$ is a \announce{$3$-transposition
group}.

The third part of the theorem is embedded in Fischer's
\cite[Lemma~(2.1.1)]{Fi71} and is also in  
\cite[Lemma~3.16]{CuHa95}.

We say that the two \trans\ groups $(G_1,D_1)$ and $(G_2,D_2)$
have the same \announce{central type} 
(usually abbreviated to \announce{type})
provided $G_1/Z(G_1)$
and $G_2/Z(G_2)$ are isomorphic as \trans\ groups.
Theorem \ref{thm-trans-prop}(c) tells us that the $3$-transposition 
properties of groups sharing a central type are
essentially the same.
In particular the two $3$-transposition groups have the same
type if and only if they have isomorphic diagrams $\diag{D_1}$ 
and $\diag{D_2}$.

\dspace{1}

A consequence of the work by Fischer \cite{Fi71} and 
later Cuypers and Hall \cite{CuHa95} is the classification 
up to isomorphism of 
all diagrams for all finite $3$-transposition
groups.\footnote{Beware: nonisomorphic groups may have the same diagram.}
In Section \ref{sec-eig} we shall give the eigenvalues 
and spectrum of
(the adjacency matrix of) each such diagram.
If $\spec{\dots,r_i,\dots}$ is the spectrum of
$\AMat{\diag{D}}$, we also say that $\spec{\dots,r_i,\dots}$ is the
spectrum of $(G,D)$ and $G$.

As $D$ is a conjugacy class of the \trans\ group $(G,D)$,
its diagram $\diag{D}$ is connected. Thus the degree
$k$ of the diagram is an eigenvalue of multiplicity
one (associated with the eigenvector $\allone$) and will
be listed first in the spectrum. The remaining eigenvalues
are restricted.

For the \trans\ group $(G,D)$ we write
$p^{\comd{h}}$, with $p \in \{2,3\}$,
for a normal $p$-subgroup $N$ with $|D \cap dN|=p^h$ for all $d \in D$.
We call this the \announce{shape} of $N$.

\begin{lem}
	\label{lem-adj-mat}
Let $(G,D)$ be a \trans\ group with normal subgroup $N$ of shape
$p^\comd{h}$ and \trans\ quotient $(H,E)$ for $H=G/N$ and $E=DN/N$.
The adjacency matrix of the diagram $\diag{G}\,(=\diag{D})$ for 
the group $G=p^\comd{h}H$ is
the matrix
$p^\comd{h}M$
of Proposition $\ref{prop-eig23}$,
where $M$ is the adjacency matrix of the diagram 
$\diag{H}\,(=\diag{E})$ for $H$.
\end{lem}

\pf Let $d,e\in D$.
If $p=2$, then the $2^h$ vertices of $\diag{dN \cap D}$ admit no edges
while if $p=3$ the subdiagram $\diag{dN\cap D}$ of size $3^h$ is complete.
In both cases, if $de$ has order $2$ then there are no edges
between $\diag{dN \cap D}$ and $\diag{eN\cap D}$, while
if $de$ has order $3$, all possible edges
between $\diag{dN \cap D}$ and $\diag{eN\cap D}$ occur. \qed

\dspace{.5}

For fixed $H$ and $p^\comd{h}$, 
there may be \trans\ groups $G$ of distinct central type 
with $N$ of type $p^\comd{h}$
and $G/N=H$, so that they have the same diagram.

Proposition \ref{prop-eig23} yields:

\begin{cor}\label{cor-upone}
  If the spectrum of $(H,E)$ is $\spec{k;\dots,r_i,\dots}$ then:
\begin{enumerate}[{\rm(a)}]
  \item the spectrum of $G=2^\comd{1}H$ is
     $\spec{2k;0, \dots, 2r_i,0, \dots}$,
     and for each $r_i \notin \{k,0\}$, the multiplicity of $2r_i$ for $G$ is
     equal to that of $r_i$ for $H$;
  \item the spectrum of $G=3^\comd{1}H$ is
     $\spec{3k+2;-1,-1, \dots, 3r_i+2,-1,-1 \dots}$,
     and for each $r_i \notin \{ k,-1\}$, the multiplicity of $3r_i+2$ for $G$ is
     equal to that of $r_i$ for $H$.
\qed
\end{enumerate}
\end{cor}

\noindent
The first part of the corollary appears in Matsuo's
original papers on vertex operator algebras 
\cite[Lemma~4.1.3]{Ma03} and \cite[\S 5]{Ma05}.

As a first example, the \trans\ group $\Sym{2}$ has diagram adjacency matrix
$M=[0]$ with unique eigenvalue $k=0$. Therefore $\Sym{3}=3^\comd{1}\Sym{2}$
has diagram adjacency matrix
\[
  \begin{bmatrix}
    0&1&1\\ 1&0&1\\ 1&1&0
  \end{bmatrix}
\]
with spectrum 
$\spec{3\cdot0+2; -1,-1}=\spec{2;-1,-1}$ and 
$\Sym{4}=2^\comd{1}\Sym{3}$ has spectrum 
\begin{eqnarray*}
&\spec{2\cdot2; 0, 2\cdot (-1), 0, 2\cdot(-1), 0}
=\spec{4;0, -2,0, -2,0}\\
&=\spec{4;[-2]^2,[0]^3} =\spec{4;[-2]^2,[0]^\star}\,.
\end{eqnarray*}
Here again we use the convention that $[t]^c$ indicates
an eigenvalue $t$ of multiplicity $c$. We also introduce the
notation $[t]^\star$ to indicate that the eigenvalue $t$
has multiplicity equal to whatever is required for the total
multiplicity to be the size $n$.

We can continue in this fashion, so that
\[
  \su{3}'=3^\comd{2}\Sym{2}=3^\comd{1}(3^\comd{1}\Sym{2})
\]
has spectrum
\[
  \spec{8;-1,-1,\,-1,-1,-1,\,-1,-1,-1}=
  \spec{8;[-1]^8}
  = \spec{8;[-1]^\star}\,
\]
while
\[
  \Weyl{D_4}=2^\comd{2}\Sym{3}=2^\comd{1}(2^\comd{1}\Sym{3})
\]
has spectrum
\[
\spec{8;0\,,\, -4,0\,,\, -4,0\,,\, 0,0\,,\, 0,0\,,\, 0,0}=
\spec{8;[-4]^2,[0]^9}=
\spec{8;[-4]^2,[0]^\star}\,.
\]

Iteration of the previous corollary yields:

\begin{cor}\label{cor-step}
	Let the $3$-transposition group $(H,E)$ have spectrum 
	\[	
	\spec{k;\dots,\eigmult{r_i}{m_i},\dots }
	\]	
	and size is
	$n_H=|E|=1 + \sum_i m_i$.
	\begin{enumerate}[$(a)$]
	\item
	A $3$-transposition group
	$G=2^\comd{h}H$, for $h \ge 1$,
	with class $D_G=2^\comd{H}E$ 
	has size
	\[
	n_G=|D_G|=2^h n_H
	\]
	and spectrum
	\[	
	\spec{2^hk;\dots,\eigmult{2^hr_i}{m_i},\dots, \eigmult{0}{\star}}
	\,.
	\]
	\item
	A $3$-transposition group
	$F=3^\comd{h}H$, for $h \ge 1$,
	with class $D_F=3^\comd{H}E$ 
	has size
	\[
	n_F=|D_F|=3^h n_H
	\]
	and spectrum
	\[	
	\spec{3^h(k+1)-1;
    \dots,\eigmult{3^h(r_i+1)-1}{m_i},\dots, \eigmult{-1}{\star}}
	\,.
	\]
	\qed
	\end{enumerate}	
\end{cor}

Note that in (a) one of the $r_i$ may be zero, in which case the expected 
tail multiplicity $(2^h-1)n_H$ should be combined with the multiplicity 
$m_i$ of $2^hr_i=0$; this explains the exponent $\star$, which indicates
a multiplicity that is whatever is needed to exhaust all eigenvalues.
Similarly, in (b) one of $r_i$ may be $-1$ and then the expected
tail multiplicity $(3^h-1)n_H$ is added to the multiplicity $m_i$ 
of $3^h(r_i+1)-1=-1$.

\section{Classifications of $\mathbf 3$-transposition groups}

Fischer's \cite{Fi71} main theorem on \trans\ groups is:

\begin{thm}\label{thm-fischer}
Let $(G,D)$ be a finite \trans\ group with no nontrival
normal solvable subgroup. Then
the group $G$ has exactly one of the central types below.
Furthermore, for each $G$ the generating class $D$ is uniquely 
determined up to an automorphism of $G$.

\name{I2.} 
$\Sym{m}$, all $m\ge 5$;

\name{I3.} $\ogpt{2m}$, $\eps=\pm$, 
all $m \ge 3$, $(m,\eps) \neq (3,+)$;

\name{I4.} $\spt{2m}$, all $m \ge 3$;

\name{I5.} $\omgp{m}{\eps}$, $\eps=\pm$, all $m \ge 6$;

\name{I6.} $\su{m}$, all $m \ge 4$;

\name{I7.}
$\Fi{22}, \Fi{23}, \Fi{24}$,
$\tri{2}, \tri{3}$.
\end{thm}

\noindent
The notation is that of \cite{CuHa95} and will be discussed
in the next section.

No example appears twice in the theorem.
Apparent omissions within it, in the first table
of the next section, and throughout the paper are explained by
the following coincidences.

\begin{lem}
\label{lem-except} 
\begin{enumerate}[$(a)$]
\item $1= \omgp{1}{+}$;
\item $\Sym{2} \iso \omgp{1}{-} \iso \omgp{2}{+} \iso \oZ_2$:
\item 
$1 \neq \oO_3(\Sym{3})$; 
$\Sym{3} \iso \ogpt[-]{2} \iso \spt{2} \iso \su{2}$;
\item 
$1 \neq \oO_2(\Sym{4})$; 
$\Sym{4} \iso \omgp{3}{-}$;
\item $\Sym{5} \iso \ogpt[-]{4}$;
\item $\Sym{6} \iso \spt{4}\iso\omgp{4}{-}$;
\item $\Sym{8} \iso \ogpt[+]{6}$;
\item $\ogpt[-]{6} \iso \omgp{5}{+}$;
\item $2 \times \su{4}\iso\omgp{5}{-}$;
\item $1 \neq \oO_2(G)$ for $G \in
\{\ogpt[+]{4},\omgp{2}{-},\omgp{3}{+},\omgp{4}{+}\}$;
\item $1 \neq \oO_3(\su{3}')$.
\end{enumerate}
\end{lem}

\pf All these can be found in \cite[\S 2]{CuHa95}. \qed


\dspace{.5}

Fischer's theorem was extended in \cite{CuHa95}. A consequence
of the main theorem of that paper is:

{\allowdisplaybreaks
\begin{thm}\label{thm-cuha}
Let $(G,D)$ be a finite \trans\ group.
Then, for integral $m$ and $h$, the group $G$ has one of the 
central types below.
Furthermore, for each $G$ the generating class $D$ is uniquely 
determined up to an automorphism of $G$.

\name{PR1.} $3^{\comd{h}}\splt\Sym{2}$, all $h \ge 1$;

\name{PR2(a).} $2^\comd{h}\splt\Sym{m}$, all $h \ge 0$, all $m\ge 4$;

\name{PR2(b).} $3^\comd{h}\splt\Sym{m}$, all $h \ge 1$, all $m\ge 4$;

\name{PR2(c).} $3^\comd{h}\splt2^\comd{1}\splt\Sym{m}$, all $h \ge 1$, 
all $m\ge 4$;

\name{PR2(d).} $4^\comd{h}\splt3^\comd{1}\splt\Sym{m}$, all $h \ge 1$, 
all $m\ge 4$;

\name{PR3.} $2^\comd{h}\splt\ogpt{2m}$, 
$\eps=\pm$, all $h \ge 0$, all $m \ge 3$, $(m,\eps) \neq (3,+)$;

\name{PR4.} $2^\comd{h}\splt\spt{2m}$, all $h \ge 0$, all $m \ge 3$;

\name{PR5.} $3^\comd{h}\,\omgp{m}{\eps}$, $\eps=\pm$,
all $h \ge 0$, all $m \ge 5$;

\name{PR6.} $4^\comd{h}\su{m}'$, all $h \ge 0$, all $m \ge 3$;

\name{PR7(a-e).} $\Fi{22}, \Fi{23}, \Fi{24}$, $\tri{2}, \tri{3}$;

\name{PR8.} $4^\comd{h}\splt(\osix)$, all $h \ge 1$;

\name{PR9.} $3^\comd{h}\splt(2 \times \spt{6})$, all $h \ge 1$;

\name{PR10.} $3^\comd{h}\splt(2\nsplt\ogpt[+]{8})$, all $h \ge 1$;

%
\name{PR11.} $3^\comd{2h}\splt(2 \times \su{5})$, all $h \ge 1$;

\name{PR12.} $3^\comd{2h}\splt4^\comd{1}\splt\su{3}'$, all $h \ge 1$.
\end{thm}
}

The notation \name{Ik} and \name{PRk} of the two theorems 
and the tables of the next section comes
from \cite{CuHa95}, where the first suggests that the groups
act \name{I}rreducibly on their natural modules, while the
second says that more general examples arise from \name{P}arabolic
subgroups of the irreducible examples---specifically their
subgroups generated by \name{R}eflections or transvections,
as appropriate.

\dspace{.5}

In the theorem
(and elsewhere) $A\splt B$ indicates a split group extension 
with normal subgroup $A$ while $A\nsplt B$ is a nonsplit group 
extension with normal subgroup $A$ and quotient $B$. The related 
notation $AB$ indicates that $A$ is normal with quotient $B$, but 
the extension may or may not be split. Extensions are left-adjusted, 
so in $A \splt B \splt C$, the normal subgroup $A \splt B$ is split 
by $C$ while $A \splt B$ has $A$ normal and split by $B$.

Neither the actual structure of the normal $p$-subgroup
nor the splitting of the extension affect the shape
of the normal subgroup and so the diagram. 
This allows us in the theorem to bundle the
exotic cases \name{PR13-19} from \cite{CuHa95}
under the corresponding generic cases \name{PR5-6}
(where both split and nonsplit group extensions may occur).

In Theorem \ref{thm-cuha}
we have rewritten shapes $2^\comd{2h}$ as $(4)^\comd{h}$
when all the nontrivial composition factors 
in the normal subgroup 
those factors
are naturally $\gf{4}$-modules for the quotient.

The only three repetitions on the list are 
$3^\comd{2}\splt\Sym{2} \iso \su{3}'$
appearing under both \name{PR1} and \name{PR6};
$\omgp{5}{+} \iso \ogpt[-]{6}$
appearing in the $h=0$ cases of both \name{PR5} and \name{PR3};
and $\omgp{5}{-} \iso 2 \times \su{4}$
appearing in the $h=0$ cases of both \name{PR5} and \name{PR6}.

Under \name{PR2(a)} the groups $2^\comd{h}\splt\Sym{3}$, for $h \ge 1$,
have the same central type as $2^\comd{h-1}\splt\Sym{4}$. Other apparent
absences are justified by Lemma \ref{lem-except}.

\section{Case analysis of spectra}
\label{sec-eig}

In Theorem \ref{thm-cuha}, each choice of parameters in
each part yields a unique diagram which admits a
$3$-transposition group (and perhaps many). In this
section we calculate the size (number of vertices) and
spectrum of the diagram in each case.
These are collected in two tables---one for the Irreducible
	examples of Theorem \ref{thm-fischer} and a second for the
	Parabolic Reflection examples of Theorem \ref{thm-cuha}. 
	The second of these essentially comes from combining the first 
	with Corollary \ref{cor-step}.

\dspace{.5}

As Fischer noted \cite[Theorem 3.3.5]{Fi71}, in each case of Theorem 
\ref{thm-fischer} (except for the triality groups $\tri{2}$ and 
$\tri{3}$), the permutation representation
of $G$ acting on $D$ by conjugation is primitive of rank $3$.
Therefore the corresponding spectrum obeys all the conditions
discussed in Section \ref{sec-rk3}. 

The redundancy of the parameter sets is of aid here.
We have $n=|D|$. As the codiagram is the commuting graph of $D$,
we also have
\[
k = |\oC_D(d)|-1
\]
for $d \in D$,
where $\oC_D(d)\backslash \{d\}$ is the  noncentral normal set (indeed
conjugacy class) of $3$-transpositions in the subgroup
$\oC_G(d)$. Similarly if $e \in \oC_D(d)\backslash \{d\}$
and $c \in D \backslash \oC_D(d)$, then
\[
\lam' =  |\oC_D(d,e)|-2
\quad\text{and}\quad
\mu'= |\oC_D(d,c)|
\]
count the $3$-transpositions of the subgroups $\oC_G(d,e)$
and $\oC_G(d,c)$.

\dspace{.3}

We have seen in Theorems \ref{thm-fischer} and \ref{thm-cuha} that
in the pair $(G,D)$ the group $G$ determines the generating class
$D$ uniquely up to an 
automorphism. Therefore we may abuse notation
by writing $\diag{G}$ for the diagram in place of $\diag{D}$.

Most of the results given here could be extracted from the
literature---for instance \cite{Hu75} and \cite{aeb}---although
the notation varies enough that translation into the form
we desire can be difficult. 
We have recalculated everything (to our own satisfaction) but
only outline the paths taken.

The first table gives the extended parameters 
$(n,k,\lam,\mu\,;\{\eigmult{r}{f},\eigmult{s}{g}\})$
of the rank $3$ (strongly regular) 
codiagrams $\codiag{G}$ and diagrams $\diag{G}$. 
Note the set notation 
for the eigenvalues and their multiplicities. This is because in some 
cases the roles of $r$ (positive eigenvalue) and $s$ (negative eigenvalue) 
may switch depending on the value of $m$. In these cases we use $d$ and $e$ 
for multiplicities in order to avoid misleading the reader.

The second table gives the size $n$ and spectrum 
$\spec{k;\dots,\eigmult{r_i}{m_i},\dots }$ of all 
diagrams $\diag{G}$. The eigenvalue in bold is the minimum eigenvalue. 
This will be of relevance in Section 7.

In Theorem \ref{thm-cuha} we have restricted parameters in
order to minimize repetition of examples. In the second
table we reverse that decision, enlarging the parameter sets
to a natural level of generality. In particular, unless
otherwise stated, $h$ can be any nonnegative integer.
\newpage
\thispagestyle{empty}
\dspace{-7}
\[
\hspace{-1cm}
\begin{array}{|c|c|}
\hline
\text{\bfseries Graph} & \text{\bfseries Extended parameters}
\strut\ (n,k,\lam,\mu\,;\{\eigmult{r}{f},\eigmult{s}{g}\})
\\
\hline\hline
\text{\name{I2}}
& 
\\ \hline
\codiag{\Sym{m}},\
m \ge 4
& 
		  \left(\binom{m}{2},\,
		  \binom{m-2}{2},\,\binom{m-4}{2},\binom{m-3}{2}\,; 
          \{
		  \eigmult{1}{m(m-3)/2},
		  \eigmult{-m+3}{m-1}
		  \}
		  \right)
\\ \hline 
\diag{\Sym{m}}
\,,\
m \ge 4
& 
\left(\binom{m}{2}\,,\,2(m-2)\,,\,m-2\,,4\,; \{
\eigmult{m-4}{m-1},\eigmult{-2}{m(m-3)/2}\}\right)
\\ \hline\hline
\text{\name{I3}}
& 
\\ \hline
\codiag{\ogpt[\eps]{2m}}
&
\bigl(
(2^{2m-1}-\eps2^{m-1},\,2^{2m-2}-1,\,
2^{2m-3}-2,\,
2^{2m-3}+\eps 2^{m-2}
\,;
\\
\eps = \pm\,,\
m \ge 2
&
\{
\eigmult{\eps2^{m-2}-1}{(2^{2m}-4)/3},
\eigmult{-\eps 2^{m-1}-1}{(2^m-\eps 1)(2^{m-1}-\eps1)/3}
\}
\bigr)
\\ \hline
\diag{\ogpt[\eps]{2m}}
& 
		\bigl(
		2^{2m-1}-\eps2^{m-1},\,
		2^{2m-2}-\eps2^{m-1},\,
		2^{2m-3}-\eps2^{m-2},\,
		2^{2m-3}-\eps2^{m-1}
		\,;
		\\
		\eps = \pm\,,\
		m \ge 2
&
		\{
		\eigmult{\eps 2^{m-1}}{(2^m-\eps 1)(2^{m-1}-\eps1)/3},
		\eigmult{-\eps2^{m-2}}{(2^{2m}-4)/3}
		\}
		\bigr)
\\ 
\hline\hline
\text{\name{I4}}
& 
\\ \hline
\codiag{\spt{2m}}
& 
		\bigl(2^{2m}-1\,,\,2^{2m-1}-2,\,2^{2m-2}-3,\,2^{2m-2}-1\,;
		\\
		m \ge 2 &
		\{
		\eigmult{2^{m-1}-1}{2^{2m-1}+2^{m-1}-1},
		\eigmult{-2^{m-1}-1}{2^{2m-1}-2^{m-1}-1}
		\}
		\bigr)
\\ \hline
\diag{\spt{2m}}
& 
	\bigl(2^{2m}-1\,,\,2^{2m-1},\,2^{2m-2},\,2^{2m-2}\,; 
	\\
	m \ge 2&
	\{
	\eigmult{2^{m-1}}{2^{2m-1}-2^{m-1}-1},\eigmult{-2^{m-1}}{2^{2m-1}+2^{m-1}-1}
	\}
	\bigr)
\\ \hline\hline
\text{\name{I5}}
& 
\\ \hline
\codiag{\omgp{m}{\eps}}
&
\bigl(
(3^{m-1}-\eps3^{(m-1)/2})/2\,,\,
(3^{m-2}+\eps3^{(m-3)/2})/2,
\\\text{odd}\ m \ge 5 
&
(3^{m-3}+\eps3^{(m-3)/2})/2\,,\,
(3^{m-3}+\eps3^{(m-3)/2})/2
\,;\\ \eps = \pm&
\{
\eigmult{3^{(m-3)/2}}{g}\,,\,
\eigmult{-3^{(m-3)/2}}{f}
\}
\bigr)
\\ \hline
\diag{\omgp{m}{\eps}}
& 
\bigl(
(3^{m-1}-\eps3^{(m-1)/2})/2\,,\,
3^{m-2}-2\eps3^{(m-3)/2}-1\,,
\\ 
\text{odd}\ m \ge 5 &
2(3^{m-3}-\eps3^{(m-3)/2}-1)\,,\,
2(3^{m-3}-\eps3^{(m-3)/2})
\,;\\ \eps = \pm&
\{
\eigmult{3^{(m-3)/2}-1}{f}\,,\,
\eigmult{-3^{(m-3)/2}-1}{g}
\}
\bigr)
\\ \hline
\text{with} &
\qquad f=(3^{m-1}-1-(\eps-1)(3^{(m-1)/2}-1))/4
\\&
\qquad
g=(3^{m-1}-1-(\eps+1)(3^{(m-1)/2}+1))/4
\\ \hline
\codiag{\omgp{m}{\eps}}
&
	\bigl(
	(3^{m-1}-\eps3^{(m-2)/2})/2\,,\,
	(3^{m-2}-\eps3^{(m-2)/2})/2\,,
	\\ \text{even}\ m \ge 6 &
	(3^{m-3}+\eps3^{(m-4)/2})/2\,,\,
	(3^{m-3}-\eps3^{(m-2)/2})/2
	\,;
	\\\eps = \pm &
	\{\,
	\eigmult{\eps3^{(m-2)/2}}{d}\,,\,\eigmult{-\eps3^{(m-4)/2}}{e}
	\,\}
	\bigr)\,,
\\ \hline
\diag{\omgp{m}{\eps}}
&
\bigl(
(3^{m-1}-\eps3^{(m-2)/2})/2\,,\,
3^{m-2}-1\,,
\\ \text{even}\ m \ge 6  &
2(3^{m-3}-1)\,,\,
2(3^{m-3}+\eps3^{(m-4)/2})
\,;
\\ \eps = \pm &
\{\,
\eigmult{-\eps3^{(m-2)/2}-1}{d}\,,\,\eigmult{\eps3^{(m-4)/2}-1}{e}
\,\}
\bigr)
\\ \hline
\text{with} &
\qquad
d=(3^{m/2}-\eps)(3^{(m-2)/2}-\eps)/8
\\&
\qquad
e=(3^m-9)/8
\\ \hline\hline
\text{\name{I6}}
& 
\\ \hline
\codiag{\su{m}}
& 
		\bigl(
		(2^{2m-1}-(-2)^{m-1}-1)/3\,,\,
		4(2^{2m-5}-(-2)^{m-3}-1)/3\,,
		\\m \ge 4 &
	  3+16(2^{2m-9}-(-2)^{m-5}-1)/3\,,\,
		(2^{2m-5}-(-2)^{m-3}-1)/3
		\,;\\&
		\{
		\eigmult{(-2)^{m-3}-1}{d},
		\eigmult{(-2)^{m-2}-1}{e}
		\}
		\bigr)
\\ \hline 
\diag{\su{m}}
& 
			(2^{2m-1}-(-2)^{m-1}-1)/3\,,\,
			2^{2m-3}\,,
			\\
			m \ge 4&
			3(2^{2m-5}) + (-2)^{m-3}\,,\,
			3(2^{2m-5})
			\,; \\ &
			\{
			\eigmult{-(-2)^{m-3}}{d},
			\eigmult{-(-2)^{m-2}}{e}
			\}
			\bigr)
\\ \hline 
\text{with}
&\qquad
d=8(2^{2m-3}-(-2)^{m-2}-1)/9
\\
&\qquad
e=4(2^{2m-3}-7(-2)^{m-3}-1)/9
\\ \hline\hline
\text{\name{I7}}
& 
\\ \hline
\codiag{\Fi{22}}
& 
	\left(3510,693,180,126\,; 
	\{
	\eigmult{63}{429},
	\eigmult{-9}{3080}
	\}
	\right)
\\ \hline
\diag{\Fi{22}}
& 
  \left(3510,2816,2248,2304\,;
  \{ 
  \eigmult{8}{3080},\eigmult{-64}{429}
  \}
  \right)
\\ \hline
\codiag{\Fi{23}}
& 
	\left(31671,3510,693,351\,;
	\{ 
	\eigmult{351}{782},\eigmult{-9}{30888}
	\}
	\right)
\\ \hline
\diag{\Fi{23}}
& 
	\left(31671,28160,25000,25344\,;
	\{ 
	\eigmult{8}{30888},\eigmult{-352}{782}
	\}
	\right)
\\ \hline
\codiag{\Fi{24}}
& 
	\left(306936,31671,3510,3240\,;
	\{ 
	\eigmult{351}{57477},
	\eigmult{-81}{249458}
	\}
	\right)
\\ \hline
\diag{\Fi{24}}
& 
	\left(306936,275264,246832,247104\,;
	\{ 
	\eigmult{80}{249458},\eigmult{-352}{57477}
	\}
	\right)
\\ \hline
\end{array}
\]
\newpage
\thispagestyle{empty}
\dspace{-7}	
\[
\hspace{-4cm}
\begin{array}{|c|c|}
\hline
\text{\bfseries Label} & \text{\bfseries Diagram}
\ \diag{G}
\\
\hline
\text{\bfseries Size} 
\ n
& \text{\bfseries Spectrum}
\ \spec{k;\dots,\eigmult{r_i}{m_i},\dots }\\
\hline\hline
\text{\name{PR1}} & \diag{3^\comd{h}\splt\Sym{2}}\\
\hline
3^h & \spec{3^h-1;{[\mathbf{-1}]^{-1+3^h}}}\\
\hline\hline
\text{\name{PR2(a)}} & \diag{2^\comd{h}\splt\Sym{m}}
\\
\hline
m \ge 4:\,
2^{h-1}m(m-1) & 
	\spec{2^{h+1}(m-2);\eigmult{2^h(m-4)}{m-1},\eigmult{0}{\star},
	\eigmult{\mathbf{-2^{h+1}}}{m(m-3)/2}}
	\\
\hline
\text{\name{PR2(b)}} & \diag{3^\comd{h}\splt\Sym{m}}
\\
\hline
m \ge 4\,, h \ge 1:\,
3^hm(m-1)/2 & \spec{3^{h}(2m-3)-1;
	\eigmult{3^h(m-3)-1}{m-1},[-1]^\star,
	\eigmult{\mathbf{-3^h-1}}{m(m-3)/2}}
	\\
\hline
\text{\name{PR2(c)}} & \diag{3^\comd{h}\splt 2^\comd{1} \splt\Sym{m}}
\\
\hline
m \ge 4\,, h \ge 1:\,
3^{h}m(m-1) & 	\spec{3^h(4m-7)-1;
	\eigmult{3^h(2m-7)-1}{m-1},
	\eigmult{3^h-1}{m(m-1)/2},\eigmult{-1}{\star}, 
	\eigmult{\mathbf{-3^{h+1}-1}}{m(m-3)/2}
		}\\
\hline
\text{\name{PR2(d)}} & \diag{4^\comd{h}\splt 3^\comd{1} \splt\Sym{m}}
\\
\hline
m \ge 4\,, h \ge 1:\
3(2^{2h-1})m(m-1) & 	\spec{4^h(6m-10);
	\eigmult{4^h(3m-10)}{m-1},\eigmult{0}{\star},\eigmult{-4^h}{m(m-1)},
	\eigmult{\mathbf{-4^{h+1}}}{m(m-3)/2}}\\
\hline\hline
\text{\name{PR3}} & \diag{2^\comd{h}\splt\ogpt[\eps]{2m}}
\\
\hline
m \ge 3\,,\ \eps=+:\,
2^{h}(2^{2m-1}-2^{m-1})
& 			
\spec{
	2^h(2^{2m-2}-2^{m-1});\,
	\eigmult{2^{h+m-1}}{(2^m-1)(2^{m-1}-1)/3},
	\eigmult{0}{\star},
	\eigmult{\mathbf{-2^{h+m-2}}}{(2^{2m}-4)/3},
	}
 \\
\hline
m \ge 2\,,\
\eps=-:\,
2^{h}(2^{2m-1}+2^{m-1})
& 			
\spec{
	2^h(2^{2m-2}+2^{m-1});\,
		\eigmult{2^{h+m-2}}{(2^{2m}-4)/3},
		\eigmult{0}{\star},
	\eigmult{\mathbf{-2^{h+m-1}}}{(2^m+1)(2^{m-1}+1)/3}
	}
\\
\hline\hline
\text{\name{PR4}} & \diag{2^\comd{h}\splt\spt{2m}}
\\
\hline
m \ge 2\,,\
2^{h}(2^{2m}-1) & 	\spec{2^{2m-1+h};\eigmult{2^{m-1+h}}{2^{2m-1}-2^{m-1}-1},
	[0]^\star,
	\eigmult{\mathbf{-2^{h+m-1}}}{2^{2m-1}+2^{m-1}-1}}
\\
\hline\hline
\text{\name{PR5}\ (and \name{PR13-16})}& 
\diag{3^\comd{h}\,\omgp{m}{\eps}}
\\
\hline
\text{odd}\ m \ge 5
& 
\spec{3^h(3^{m-2}-2\eps3^{(m-3)/2})-1\,;\,
        \eigmult{3^{(m-3)/2+h}-1}{f}\,,\,\eigmult{-1}{\star},
        \eigmult{\mathbf{-3^{(m-3)/2+h}-1}}{g}
}\\
\hline
\eps=+:\,
3^h(3^{m-1}-3^{(m-1)/2})/2
&	\text{for}\ 
f=(3^{m-1}-1)/4
\ \text{and}\
g=(3^{m-1}-1-2(3^{(m-1)/2}+1))/4
\\
\hline
\eps=-:\,
3^h(3^{m-1}+3^{(m-1)/2})/2
&	\text{for}\
f=(3^{m-1}-1+2(3^{(m-1)/2}-1))/4
\ \text{and}\
g=(3^{m-1}-1)/4
\\ 
\hline
\text{even}\ m \ge 6
\,,\ \eps = +:\,
& 
\spec{3^{m-2+h}-1\,
        ;\,
        \eigmult{3^{(m-4)/2+h}-1}{f}
        \,,\,
        \eigmult{-1}{\star},
        \eigmult{\mathbf{-3^{(m-2)/2+h}-1}}{g}}
\\
\hline
3^h(3^{m-1}-3^{(m-2)/2})/2
& 
\text{for}\
f=(3^m-9)/8
\ \text{and}\
g=(3^{m/2}-1)(3^{(m-2)/2}-1)/8
\\
\hline
\text{even}\ m \ge 6
\,,\ \eps = -:\,
& 
\spec{3^{m-2+h}-1\,
        ;\,\eigmult{3^{(m-2)/2+h}-1}{f}\,,
        \eigmult{-1}{\star},
        \eigmult{\mathbf{-3^{(m-4)/2+h}-1}}{g}}
\\
\hline
3^h(3^{m-1}+3^{(m-2)/2})/2
& 
\text{for}\
f=(3^{m/2}+1)(3^{(m-2)/2}+1)/8
\ \text{and}\
g=(3^m-9)/8
\\
\hline\hline
\text{\name{PR6}\ (and \name{PR17-19})}
& 
\diag{4^\comd{h}\su{m}'}
\,,\quad 
	m \ge 3
\\ \hline
\text{even}\ m \ge 4:\,
4^{h}(2^{2m-1}-1 + 2^{m-1})/3
& 
\spec{2^{2h+2m-3};
	\eigmult{2^{2h+m-3}}{f},
	\eigmult{0}{\star},
	\eigmult{\mathbf{-2^{2h+m-2}}}{g}
	}\\
\hline
&	\text{for}\ f=8(2^{2m-3}-1-2^{m-2})/9
\ \text{and}\
g=4(2^{2m-3}-1+7(2^{m-3}))/9
\\ 
\hline
\text{odd}\ m \ge 3:\
4^{h}(2^{2m-1}-1 - 2^{m-1})/3
& 
\spec{2^{2h+2m-3};
	\eigmult{2^{2h+m-2}}{f},
	\eigmult{0}{\star},
	\eigmult{\mathbf{-2^{2h+m-3}}}{g}}\\
\hline
&	\text{for}\ 
f=4(2^{2m-3}-1-7(2^{m-3}))/9
\ \text{and}\
g=8(2^{2m-3}-1+2^{m-2})/9
\\ \hline\hline
\text{\name{PR7(a)}} 
& 
\diag{\Fi{22}}
\\ \hline
3510
& 
	\spec{2816\,; 
		\eigmult{8}{3080},\eigmult{\mathbf{-64}}{429}}
\\ \hline
\text{\name{PR7(b)}} 
& 
\diag{\Fi{23}}
\\ \hline
31671
& 
	\spec{
		28160\,;
		\eigmult{8}{30888},\eigmult{\mathbf{-352}}{782}}
\\ \hline
\text{\name{PR7(c)}} 
& 
\diag{\Fi{24}}
\\ \hline
306936
& 
	\spec{275264\,; 
		\eigmult{80}{249458},\eigmult{\mathbf{-352}}{57477}}
\\ \hline
\text{\name{PR7(d)}} 
& 
\diag{\tri{2}}
\\ \hline
360 
& 
\spec{296;\eigmult{8}{105},\eigmult{-4}{252},
	\eigmult{\mathbf{-64}}{2}}
\\ \hline
\text{\name{PR7(e)}} 
& 
\diag{\tri{3}}
\\ \hline
3240 
&
\spec{2888;[8]^{2457},[-28]^{780},[\mathbf{-352}]^2}
\\ \hline\hline
\text{\name{PR8}} 
& 
\diag{4^\comd{h}\splt(\osix)}
\\ \hline
h \ge 1:\,
126(4^h)& 
\spec{5(4^{h+2});\eigmult{2^{2h+3}}{35},\eigmult{0}{\star},
	\eigmult{\mathbf{-4^{h+1}}}{90}
	}
\\ \hline\hline
\text{\name{PR9}} & 
\diag{3^\comd{h}\splt(2 \times \spt{6}}
\\ \hline
h \ge 1:\,
63(3^{h})
& 
	\spec{11(3^{h+1})-1;\eigmult{5(3^h)-1}{27},
		\eigmult{-1}{\star},
		\eigmult{\mathbf{-3^{h+1}-1}}{35}
		}
\\ \hline\hline
\text{\name{PR10}} & 
\diag{3^\comd{h}\splt(2\nsplt\ogpt[+]{8})}
\\ \hline
h \ge 1:\,
120(3^{h})
& 
	\spec{19(3^{h+1})-1;\eigmult{3^{h+2}-1}{35},
		\eigmult{-1}{\star},
		\eigmult{\mathbf{-3^{h+1}-1}}{84}}
\\ \hline\hline
\text{\name{PR11}} & 
\diag{3^\comd{2h}\splt(2 \times \su{5})}
\\ \hline
h \ge 1:\,
165(3^{2h})
& 
	\spec{
		43(3^{2h+1})-1; \eigmult{3^{2h+2}-1}{44},
		\eigmult{-1}{\star},
		\eigmult{\mathbf{-3^{2h+1}-1}}{120}
		}
\\ \hline\hline
\text{\name{PR12}} & 
\diag{3^\comd{2h}\splt 4^\comd{1}\splt\su{3}'}
\\ \hline
h \ge 1:\,
36(3^{2h})
& 
\spec{
	11(3^{2h+1})-1;[3^{2h}-1]^{27},
	\eigmult{-1}{\star},
	\eigmult{\mathbf{-3^{2h+1}-1}}{8}
	}
\\ \hline

\end{array}
\]

\newpage
\subsection{Moufang case}

This is the situation in which the diagram $\diag{D}$
is a complete graph. That is, there are no $D$-subgroups
of $G$ isomorphic to $\Sym{4}$. The terminology comes
from a connection with commutative Moufang loops of
exponent $3$; see \cite{CuHa95}.

For $h \ge 0$, let $N_h$ be an elementary abelian $3$-group 
of order $3^h$. Further let $d$ be an element of order
$2$ that acts on $N_h$ as inversion.
Then for $G_h=N_h \splt \langle d \rangle$
and $D_h=d^{N_h}$, the pair $(G_h,D_h)$ is a $3$-transposition
group $3^\comd{h}\splt\Sym{2}$ of Moufang type \name{PR1}.
Conversely, every finite $3$-transposition group with complete
diagram arises as $N\splt\Sym{2}$ for some normal
$3$-subgroup $N$. (Appropriate $N$ exist with arbitrarily
large nilpotence class.)

By Corollary \ref{cor-step}:

\begin{thm}\label{thm-PR1-case}
	 \name{PR1}:
	 the diagram $\diag{3^\comd{h}\splt\Sym{2}}$ 
	 for $h \ge 0$ has
	 size 
	 \[
	 n=3^h
	 \] 
	 and spectrum
	 \[
	 \spec{3^h-1;-1,-1,\,\dots,\,-1,-1,-1}
	 =
	 \spec{3^h-1;[-1]^{-1+3^h}}
	 \,.
	 \]\qed
\end{thm}

The fundamental $3$-transposition groups $\bbZ_2$ and
$\Sym{3}$ occur here as $h=0$ and $h=1$.

\subsection{Symmetric cases}

\begin{thm}\label{thm-I2}
	\begin{enumerate}[$(a)$]
		\item
		For $m\ge 4$ the codiagram 
		$\codiag{\Sym{m}}$ 
		has extended parameters
		\[
		  \left(\binom{m}{2}\,,\,
		  \binom{m-2}{2}\,,\,\binom{m-4}{2}\,,\binom{m-3}{2}\,; 
		  \eigmult{1}{m(m-3)/2},
		  		  \eigmult{-m+3}{m-1}
		  \right)
		\,.
		\]
		\item
		For $m\ge 4$ the diagram 
		$\diag{\Sym{m}}$ 
		has extended parameters
		\[
		\left(m(m-1)/2\,,\,2(m-2)\,,\,m-2\,,4\,; 
		\eigmult{m-4}{m-1},\eigmult{-2}{m(m-3)/2}\right)
		\]
		and spectrum
		\[
		\spec{2(m-2);\eigmult{m-4}{m-1},\eigmult{-2}{m(m-3)/2}}
		\,.
		\]
	\end{enumerate}
\end{thm}

\pf This is well known, but it is also easy to calculate the
basic parameters of $\codiag{\Sym{m}}$ using
$3$-transposition
properties:
\begin{enumerate}[(i)]
	\item $n =\binom{m}{2}$: the $3$-transposition class 
	is $D=(1,2)^{\Sym{m}}$.
	\item $k'=\binom{m-2}{2}$: $\oC_{\Sym{m}}((1,2))$ has type $\Sym{m-2}$.
	\item $\lam'=\binom{m-4}{2}$:
	$\oC_{\Sym{m}}((1,2),(3,4))$ has type $\Sym{m-4}$.
	\item $\mu'=\binom{m-3}{2}$: $\oC_{\Sym{m}}((1,2),(2,3))$ has type $\Sym{m-3}$.
	\qed
\end{enumerate}

Corollary \ref{cor-step} gives directly:

\begin{prop}\label{prop-PR2a-case}
	\name{PR2(a)}:
	the diagram 
	$\diag{2^\comd{h}\splt\Sym{m}}$
		with $m \ge 4$ and $h \ge 0$
	has size
	\[
	n=2^{h-1}m(m-1)
	\]
	and spectrum
	\[
	\spec{2^{h+1}(m-2);\eigmult{2^h(m-4)}{m-1},\eigmult{-2^{h+1}}{m(m-3)/2}
		,\eigmult{0}{\star}}
	\,.
	\]\qed
\end{prop}

\begin{prop}\label{prop-PR2b-case}
	\name{PR2(b):}
	the diagram 
	$\diag{3^\comd{h}\splt\Sym{m}}$, 
	with $m \ge 4$ and $h \ge 0$,
	has size
	\[
	n=3^hm(m-1)/2
	\]
	and spectrum
	\[
	\spec{3^{h}(2m-3)-1;
		\eigmult{3^h(m-3)-1}{m-1},\eigmult{-3^h-1}{m(m-3)/2}
		,\eigmult{-1}{\star}}
	\,.
	\] \qed
\end{prop}

\begin{prop}\label{prop-PR2c-case}
	\name{PR2(c):}
	the diagram 
	$\diag{3^\comd{h}\splt 2^\comd{1} \splt\Sym{m}}$,
	with $m \ge 4$ and $h \ge 0$,
	has size
	\[
	n=3^{h}m(m-1)
	\]
	and spectrum
		\begin{multline*}
	\spec{-1+3^h(4m-7);
		\eigmult{3^h(2m-7)-1}{m-1},\eigmult{-3^{h+1}-1}{m(m-3)/2}
		,\\ \eigmult{3^h-1}{m(m-1)/2},\eigmult{-1}{\star}}
	\,.
		\end{multline*}
\end{prop}
	
\pf
Apply Corollary \ref{cor-step} to the
diagram $\diag{2^\comd{1} \splt\Sym{m}}$,
which has size
	\[
	n=m(m-1)
	\]
	and spectrum
	\[
	\spec{4(m-2);\eigmult{2(m-4)}{m-1},\eigmult{-4}{m(m-3)/2}
		,\eigmult{0}{m(m-1)/2}}
	\,.
	\]
\qed	

\begin{prop}\label{prop-PR2d-case}
	\name{PR2(d):}
	the diagram 
	$\diag{2^\comd{2h}\splt 3^\comd{1} \splt\Sym{m}}$,
	with $m \ge 4$ and $h \ge 0$,
	has size
	\[
	n=3(2^{2h-1})m(m-1)
	\]
	and spectrum
	\[
	\spec{4^h(6m-10);
		\eigmult{4^h(3m-10)}{m-1},\eigmult{-4^{h+1}}{m(m-3)/2},
		\eigmult{-4^h}{m(m-1)},\eigmult{0}{\star}}
	\,.
	\] 
\end{prop}

\pf
Apply Corollary \ref{cor-step} to the
diagram $\diag{3^\comd{1} \splt\Sym{m}}$,
which has size
	\[
	n=3m(m-1)/2
	\]
	and spectrum
	\[
	\spec{6m-10;
		\eigmult{3m-10}{m-1},\eigmult{-4}{m(m-3)/2}
		,\eigmult{-1}{m(m-1)}}
	\,.
	\]
\qed

\subsection{Polar space cases}
\label{sec-polar}

For us a finite polar space graph $\codiag{X}$ has as
vertex set $X$ the isotropic\footnote{more generally, singular}
$1$-spaces for a nondegenerate reflexive sesquilinear
form $f_i$ on a finite space $V_i=\gf{q}^i$ with edges given by
perpendicularity. 
In our context, the form $f_i$ is either symplectic over $\gf{2}$ 
    or hermitian over $\gf{4}$.
By Witt's Theorem, the 
corresponding isometry group acts with rank $3$ (or less). 
There are exactly two types of $2$-spaces
spanned by isotropic vectors---totally isotropic $2$-spaces
with $q+1$ 
pairwise perpendicular
isotropic $1$-subspaces and hyperbolic $2$-spaces
with \ssh{ $s$ }
pairwise nonperpendicular
isotropic $1$-subspaces.
The hyperbolic $2$-spaces are precisely the nondegenerate
$2$-spaces containing an isotropic $1$-space.
A vertex is either adjacent to all those of a given
totally isotropic $2$-space or exactly one.

Let $s_i=1+k_i+l_i$ be the number of isotropic $1$-spaces in $V_i$.
(So $s=s_2$.)
The decomposition $V_i = V_2 \perp V_{i-2}$ can be used to
calculate the parameters. This yields recursions for the degree of $\codiag{X}$
\begin{equation}
	\label{eqn-l-recur}
k'_i=l_i= q s_{i-2}
\end{equation}
and its codegree
\begin{equation}
	\label{eqn-k-recur}
l'_i=k_i = (s_2-1)q^{i-2},
\end{equation}
hence
\begin{equation}
	\label{eqn-s-recur}
s_i = 1 + (s_2-1)q^{i-2} + q s_{i-2}\,.
\end{equation}
Here we initialize with $s_1=0$
(as nondegenerate $1$-spaces contain no
isotropic vectors), 
but $s_2$ will depend upon the type
of form under consideration.
A further consequence of the decomposition is
\begin{equation}
		\label{eqn-mu-recur}
\mu_i'= s_{i-2}\,.
\end{equation}
Therefore we have the three parameters 
$s_i$, $k_i'$, and $\mu_i'$,
from which it is (at least in principal) easy to 
calculate all parameters
of $\codiag{X}$
and $\diag{X}$ using the identities of
Section \ref{sec-rk3}. 
The additional identity
\begin{equation}
    \label{eqn-lam-recur}
    \lam'_i= (q-1)+ q^2s_{i-4}\,,
\end{equation}
can be seen within $\langle x,y\rangle^\perp$, where $\langle x\rangle$ 
and $\langle y\rangle$ are distinct perpendicular isotropic $1$-spaces. 
This is because $\langle x,y\rangle^\perp/\langle x,y\rangle\cong V_{i-4}$.

\subsubsection{Symplectic over $\gf{2}$}

The nondegenerate form $f=f_{2m}$ above 
is symplectic on $V_{2m}$ if it is bilinear with
all $1$-spaces isotropic:
$f(x,x)=0$ for all $x \in V_{2m}$. Its
polar graph is denoted $\codiag{\spt[q]{2m}}$.

In the special case of symplectic polar spaces over $\gf{2}$
the corresponding
transvection isometries $D$ form a class of $3$-transpositions
in the full isometry group $G=\spt{2m}$ with the codiagram
$\codiag{D} = \codiag{X}=\codiag{\spt{2m}}$.
In this case $s_i=n_i$, $i=2m$, $q=2$, and $n_2=1+2=3$. 

\begin{thm}\label{thm-I4}
	\begin{enumerate}[$(a)$]
		\item
		For $m\ge 2$ the codiagram $\codiag{\spt{2m}}$ has extended
		parameters
		\begin{multline*}
		\bigl(2^{2m}-1\,,\,2^{2m-1}-2,\,2^{2m-2}-3,\,2^{2m-2}-1\,; \\
		\eigmult{2^{m-1}-1}{2^{2m-1}+2^{m-1}-1},
		\eigmult{-2^{m-1}-1}{2^{2m-1}-2^{m-1}-1}
		\bigr)\,.
		\end{multline*}
		\item
		For $m\ge 2$ the diagram $\diag{\spt{2m}}$ has extended
		parameters
		\[
		\bigl(2^{2m}-1\,,\,2^{2m-1},\,2^{2m-2},\,2^{2m-2}\,; 
		\eigmult{2^{m-1}}{2^{2m-1}-2^{m-1}-1},\eigmult{-2^{m-1}}{2^{2m-1}+2^{m-1}-1}
		\bigr)
		\]
		and spectrum
		\[
			\spec{2^{2m-1}
			;\,\eigmult{2^{m-1}}{2^{2m-1}-2^{m-1}-1},\eigmult{-2^{m-1}}
			{2^{2m-1}+2^{m-1}-1}}
		\,.
		\]
	\end{enumerate}
\end{thm}

\pf
\begin{enumerate}[(i)]
	\item
	$n = n_{2m}= 2^{2m}-1$: all $1$-spaces are isotropic. 
	\item
	by \eqref{eqn-l-recur} 
	$k=k_{2m}=l'_{2m}= (n_2-1)q^{i-2}=2^{2m-1}$.
	\item
	by \eqref{eqn-mu-recur} 
	$\mu'=\mu'_{2m}= n_{2m-2}=2^{2m-2}-1$.
	\item 
	by \eqref{eqn-lam-recur} $\lam'=(q-1)+q^2n_{2m-4}
	=1+4(2^{2m-4}-1)=2^{2m-2}-3$.
	\qed
\end{enumerate}

Corollary \ref{cor-step} gives immediately:

\begin{prop}\label{prop-PR4-case}
	\name{PR4:}
	the diagram 
	$\diag{2^\comd{h}\splt\spt{2m}}$
	with $m \ge 2$
	and $h \ge 0$
	has size
	\[
	n=2^{h}(2^{2m}-1)
	\]
	and spectrum
	\[
	\spec{2^{2m-1+h};\eigmult{2^{m-1+h}}{2^{2m-1}-2^{m-1}-1},
		\eigmult{-2^{m-1+h}}{2^{2m-1}+2^{m-1}-1}
		,[0]^\star}
	\,.
	\]\qed
\end{prop}

\begin{prop}\label{prop-PR9-case}
	\name{PR9:}
	the diagram 
	$\diag{3^\comd{h}\splt(2 \times \spt{6}}$
	with $h \ge 0$
	has size
	\[
	n=63(3^{h})
	\]
	and spectrum
	\[
	\spec{11(3^{h+1})-1;\eigmult{5(3^h)-1}{27},
		\eigmult{-3^{h+1}-1}{35}
		,[-1]^\star}
	\,.
	\]
\end{prop}

\pf
Apply Corollary \ref{cor-step} to the
diagram $\diag{\spt{6}} = \diag{2 \times \spt{6}}$,
which has extended parameters
		\[
		\bigl(63\,,\,32,\,16,\,16\,; 
		\eigmult{4}{27},\eigmult{-4}{35}
		\bigr) \,.
		\]
\qed

\subsubsection{Unitary over $\gf{4}$}

For finite unitary polar graphs we must have
$q=t^2$ for some prime power $t$.
The nondegenerate form $f=f_{m}$ is  
hermitian (or unitary)
on $V_{m}$ if it is biadditive with
\[
f(ax,by)=af(x,y)b^t
\]
and
\[
f(x,y)=f(y,x)^t
\]
for all $x,y \in V_{m}$ and $a,b \in \gf{q}$.
Its polar graph is denoted $\codiag{\su[t]{2m}}$.

In the special case of 
unitary polar spaces over $\gf{4}$ 
the corresponding
transvection isometries $D$ form a class of $3$-transpositions
in the isometry group $G=\su{m}$ with the codiagram
$\codiag{D} = \codiag{X}=\codiag{\su{m}}$.
In this case $s_i=n_i$, $i=m$, $q=4$, $t=2$, and $n_2=1+2=3$. 

\begin{thm}
		For $m \ge 3$ set
		\[
		d=8(2^{2m-3}-1-(-2)^{m-2})/9
		\]
		and
		\[
		e=4(2^{2m-3}-1-7(-2)^{m-3})/9\,.
		\]
	\label{thm-I6}
		\begin{enumerate}[$(a)$]
			\item
			For $m\ge 3$ the codiagram $\codiag{\su{m}}$ has extended
			parameters
		\begin{multline*}
			\bigl(
			(2^{2m-1}-1 - (-2)^{m-1})/3\,,\,
			2^2(2^{2m-5}-1 - (-2)^{m-3})/3\,,
			\\
			\lam'=
			3+16(2^{2m-9}-(-2)^{m-5}-1)/3\,,\,
			\mu'=(2^{2m-5}-1-(-2)^{m-3})/3
			\,;\\
			\{\eigmult{r'}{g},\eigmult{s'}{f}\}=
			\{
			\eigmult{(-2)^{m-3}-1}{d},
			\eigmult{(-2)^{m-2}-1}{e}
			\}
			\bigr)\,,
		\end{multline*}
			\item
		For $m\ge 3$ the diagram $\diag{\su{m}}$ has extended
		parameters
		\begin{multline*}
			\bigl(
			(2^{2m-1}-1 - (-2)^{m-1})/3\,,\,
			2^{2m-3}\,,\\
			\lam=3(2^{2m-5}) + (-2)^{m-3}\,,\,
			\mu=3(2^{2m-5})
			\,; \\
			\{\eigmult{r}{f},\eigmult{s}{g}\}=
			\{
			\eigmult{-(-2)^{m-3}}{d},
			\eigmult{-(-2)^{m-2}}{e}
			\}
			\bigr)\,,
		\end{multline*}
				and spectrum
			\[
			\spec{2^{2m-3};
				\eigmult{-(-2)^{m-3}}{d},
				\eigmult{-(-2)^{m-2}}{e}
				}\,.
			\]
		\end{enumerate}
\end{thm}

\pf

\begin{enumerate}[(i)]
	\item
	$n = n_m= (2^{2m-1}-1 - (-2)^{m-1})/3$: 
	the recursion of \eqref{eqn-s-recur}
	\[
	n_i = 1 + (n_2-1)q^{i-2} + q n_{i-2}
	\]
	is initialized by
	\[
	n_1=0 = (2-1-1)/3= (2^{2-1}-1 - (-2)^{0})/3
	\]
	and
	\[
	n_2=3=(8-1+2)/3
	=(2^{4-1}-1 - (-2)^{2-1})/3
	\,.	
	\]
	\item
	 $k=l'=l_m'=(n_2-1)q^{m-2} = (3-1) 4^{m-2}$.
	\item
	$\mu' = n_{m-2} = (2^{2m-5}-1-(-2)^{m-3})/3$.
	\qed
\end{enumerate}

\begin{prop}\label{prop-PR6-case}
	\name{PR6}:
	the diagram 
	$\diag{4^\comd{h}\su{m}'}$ for
	all $h \ge 0$ and 
	all $m \ge 3$
	has size
	\[
	n=4^{h}(2^{2m-1}-1 - (-2)^{m-1})/3
	\]
	and spectrum
	\[
	\spec{2^{2m-3+2h};
	\eigmult{-(-2)^{m-3+2h}}{d},
	\eigmult{-(-2)^{m-2+2h}}{e}
		,[0]^\star}
	\]
	where
	\[
	d=8(2^{2m-3}-1-(-2)^{m-2})/9
	\]
	and
	\[
	e=4(2^{2m-3}-1-7(-2)^{m-3})/9\,.
	\]
	\qed
\end{prop}

\begin{prop}\label{prop-P11-case}
	\name{PR11:} the diagram 
	$\diag{3^\comd{2h}\splt(2 \times \su{5})}$,
	for $h \ge 0$,
		has size 
		\[
		n=165(3^{2h})
		\] 
		and spectrum
		\[
		\spec{129(3^{2h})-1;
			\eigmult{3^{2h+2}-1}{44},\eigmult{-3^{2h+1}-1}{120},
			\eigmult{-1}{\star}}\,.
		\]\qed
\end{prop}

\pf
Apply Corollary \ref{cor-step} to the
diagram $\diag{\su{5}} = \diag{2 \times \su{5}}$,
which has size
$165$
and spectrum
		\[
		\spec{128;\eigmult{8}{44},
			\eigmult{-4}{120}
		}\,.
		\]
\qed

\begin{prop}\label{prop-P12-case}
	\name{PR12:} the diagram 
    $\diag{3^\comd{2h}\splt 4^\comd{1}\splt\su{3}'}$,
	for $h \ge 0$,
	has size 
	\[
	36(3^{2h})
	\] 
	and spectrum
	\[
	\spec{33(3^{2h})-1;\eigmult{3^{2h}-1}{27},
		\eigmult{-3^{2h+1}-1}{8},
		\eigmult{-1}{\star}}\,.
	\]
\end{prop}

\pf
Apply Corollary \ref{cor-step} to the
diagram 
$\diag{4^\comd{1}\splt\su{3}'}$,
which has size 
\[
n=36
\] 
and spectrum
\[
\spec{32;\eigmult{0}{27},
	\eigmult{-4}{8}}
\,.
\]
\qed

\subsection{Nonsingular orthogonal cases over $\gf{2}$}

Let $V_{2m}=\gf{2}^{2m}$ admit the nondegenerate symplectic
form $f=f_{2m}$. An associated quadratic $q_{2m}^\eps=q$ is a map 
$q \colon V_{2m}
\longrightarrow \gf{2}$ such that 
\[
f(x,y)=q(x+y)+q(x)+q(y)
\]
for all $x,y \in V_{2m}$. The vectors $x$ with $q(x)=0$
are singular and those with $q(x) =1$ are nonsingular.
Each of the two types of symplectic $2$-spaces resolves
into two types of orthogonal $2$-spaces. A totally isotropic
$2$-space is either totally singular ($q$ is identically $0$)
or is defective---it has exactly two nonsingular vectors.
A symplectic hyperbolic $2$-space is either orthogonal
hyperbolic---a unique nonsingular vector---or is 
asingular---its only singular vector is $0$. 
Thus the isometry type of a $2$-space is uniquely
determined by the number of nonsingular vectors it 
contains---respectively $0$, $2$, $1$, $3$.

Up to isometry, the form $q$
has one of two types denoted by the Witt sign $\eps$, 
equal to $+=+1$ or $-=-1$ depending
upon whether maximal totally singular spaces have dimension $m$ or $m-1$.
The corresponding diagram $\diag{\ogpt[\eps]{2m}}$ has as vertices
the nonsingular $1$-spaces $\langle x \rangle \in V_{2m}^\eps$ 
with two adjacent when
not perpendicular. That is,  
$\diag{\ogpt[\eps]{2m}}$ 
    is the subgraph of $\diag{\spt{2m}}$ induced
on the set of
$1$-spaces that are nonsingular for $q$, and correspondingly
for $\codiag{\ogpt[\eps]{2m}}$.

The symplectic transvections centered at nonsingular vectors
form a generating conjugacy class\footnote{More accurately, $D$
	is a normal set in the orthogonal group. But the only case in 
	which it is not a generating conjugacy
	class is $\ogpt[+]{4}$ where these transvections generate
	a proper normal subgroup $\Sym{3} \times \Sym{3}$.}
$D$ of $3$-transpositions in the corresponding
orthogonal group $\ogpt[\eps]{2m}$.
\begin{thm}\label{thm-I3}
		\begin{enumerate}[$(a)$]
			\item
			For $m\ge 1$ the codiagram $\codiag{\ogpt[\eps]{2m}}$ has extended
			parameters
			\begin{multline*}
				\bigl(2^{2m-1}-\eps2^{m-1},\,2^{2m-2}-1\,,\,
				2^{2m-3}-2\,,\,
				2^{2m-3}+\eps 2^{m-2};\\
				\{
				\eigmult{\eps2^{m-2}-1}{(2^{2m}-4)/3},
				\eigmult{-\eps 2^{m-1}-1}{(2^m-\eps 1)(2^{m-1}-\eps1)/3}
				\}
				\bigr)\,.
			\end{multline*}
			\item
			For $m\ge 1$ the diagram $\diag{\ogpt[\eps]{2m}}$ has extended
			parameters
			\begin{multline*}
				\bigl(
				2^{2m-1}-\eps2^{m-1},\,
				2^{2m-2}-\eps2^{m-1},\,
				2^{2m-3}-\eps2^{m-2},\,
				2^{2m-3}-\eps2^{m-1}
				\,;\\
				\{
				\eigmult{\eps 2^{m-1}}{(2^m-\eps 1)(2^{m-1}-\eps1)/3},
				\eigmult{-\eps2^{m-2}}{(2^{2m}-4)/3}
				\}
				\bigr)
				\,,
			\end{multline*}
			and spectrum
			\[
			\spec{
			2^{2m-2}-\eps2^{m-1};\,
			\eigmult{\eps 2^{m-1}}{(2^m-\eps 1)(2^{m-1}-\eps1)/3},
			\eigmult{-\eps2^{m-2}}{(2^{2m}-4)/3}}
			\,.
			\]
		\end{enumerate}
\end{thm}

\pf
\begin{enumerate}[$(i)$]
	\item
	$n_{2m}^\eps = 2^{2m-1}-\eps2^{m-1}$: initialize with
	$n_2^+=1$ and $n_2^-=3$. 
	As $V_2^+$ contains a unique nonsingular vector,
	the decomposition $V^\eps_{2m}= V_2^+
	\perp V^\eps_{2m-2}$ gives the recursion
	\[
	n_{2m}^\eps = 3 n_{2m-2}^\eps + (2^{2m-2}-n_{2m-2}^\eps)
	 = 2^{2m-2}+2n_{2m-2}^\eps\,,
	\]
	and the result follows.
	\item
	$k'=2^{2m-2}-1$: consider $\codiag{\ogpt[\eps]{2m}}$ as
	an induced subgraph of $\codiag{\spt{2m}}$. 
	For the nonsingular
	vector $x$, a $2$-space containing $x$ and in $x^\perp$ must
	be defective---of its two $1$-spaces not containing $x$, one
	is singular and one is nonsingular. Therefore 
	\[
	k'_{\codiag{\ogpt[\eps]{2m}}} 
    = \frac{1}{2}k'_{\codiag{\spt{2m}}}
	=\frac{1}{2}(2^{2m-1}-2)
	\,.
	\]
	\item
	$\lam'=2^{2m-3}-2$: again consider $\codiag{\ogpt[\eps]{2m}}$ as
	an induced subgraph of $\codiag{\spt{2m}}$. 
	In calculating $\lam'$ for the symplectic case, and more generally
	for the polar cases over $\mathbb F_q$ as in equation \ref{eqn-lam-recur}, 
	we found
	\[
	\lam'_{\codiag{\spt{2m}}}=(q-1)+q^2n_{2m-4}
	=1+4(2^{2m-4}-1)=2^{2m-2}-3\,,
	\]
	counting the $q-1$ remaining
	isotropic $1$-spaces of the totally isotropic $2$-space $\langle x,y\rangle$
	plus the $q^2$ additional $1$-spaces of each $3$-space 
	in $\langle x,y\rangle^\perp$ on $\langle x,y\rangle$, these
	enumerated by the $1$-spaces of $\langle x,y\rangle^\perp/\langle x,y\rangle$ 
	of dimension $4$ less. 
	Here we must restrict ourselves to nonsingular vectors, so the isotropic spaces 
	through $x$ and $y$ are defective. The only singular vectors in $\langle x,y \rangle$
	are $x$ and $y$, and each $3$-space in $\langle x,y\rangle^\perp$
	on $\langle x,y\rangle$ has exactly two additional nonsingular vectors.
	Therefore the count becomes
%
	 \[
	 	\lam'_{{\codiag{\ogpt[\eps]{2m}}}}
	 	=0+2(2^{2m-4}-1)=2^{2m-3}-2\,.
	 \]
 	 \item $\mu'=2^{2m-3}+\eps 2^{m-2}$: the 
 	 decomposition $V^\eps_{2m}= V_2^-	\perp V^{-\eps}_{2m-2}$ gives
 	 $\mu'= n_{2m-2}^{-\eps}$.
	 \qed
\end{enumerate}

\begin{prop}\label{prop-P3-case}
	\name{PR3:} the diagram 
	$\diag{2^\comd{h}\splt\ogpt{2m}}$
	for $h \ge 0$, $\eps = \pm$, and $m \ge 1$ has size 
	\[
	n=2^{h}(2^{2m-1}-\eps2^{m-1})
	\] 
	and spectrum
	\[
		\spec{
			2^h(2^{2m-2}-\eps2^{m-1});\,
			\eigmult{\eps 2^{h+m-1}}{(2^m-\eps 1)(2^{m-1}-\eps1)/3},
			\eigmult{-\eps2^{h+m-2}}{(2^{2m}-4)/3},
			\eigmult{0}{\star}}
		\,.
	\]\qed
\end{prop}

\begin{prop}\label{prop-P10-case}
	\name{PR10:} the diagram 
	$\diag{3^\comd{h}\splt(2\nsplt\ogpt[+]{8})}$
	for $h \ge 0$ has size 
	\[
	n=120(3^{h})
	\] 
	and spectrum
	\[
	\spec{57(3^{h})-1;\eigmult{3^{h+2}-1}{35},
		\eigmult{-3^{h+1}-1}{84},\eigmult{-1}{\star}}\,.
	\]
\end{prop}

\pf
Apply Corollary \ref{cor-step} to the
diagram $\diag{\ogpt[+]{8}} = \diag{2\nsplt\ogpt[+]{8}}$,
which has extended parameters
\[
\left(120,56,28,24\,; 
\eigmult{8}{35},\eigmult{-4}{84}
\right)\,.
\]
\qed

\subsection{Nonsingular orthogonal cases over $\gf{3}$}

Let $V=V_{m}=\gf{3}^{m}$ admit the nondegenerate symmetric bilinear
(that is, orthogonal) form $f$. The diagonal of $f$ yields the
quadratic form $q\colon V \longrightarrow \gf{3}$ given by
\[
q(x)=f(x,x)\,.
\] Conversely, $f$ can be reconstructed from $q$
via 
\[
f(x,y)= -q(x+y)+q(x)+q(y)\,.
\] 
In this context, the
isotropic vectors (those $x$ with $q(x)=f(x,x)=0$) 
are called singular. 
As described at the beginning of Section \ref{sec-polar},
the singular ($=$ isotropic) $1$-spaces form the
vertex set of a polar space graph which is strongly regular
and indeed rank $3$ (by Witt's Theorem). There are two types
of nonsingular $1$-spaces $\langle x\rangle$; those
with $q(x)=1$ are called $+$-spaces and $x$ is a
$+$-vector; those with $q(x)=-1$ are $-$-spaces 
and $x$ is a $-$-vector.

There are two parameters of interest for the nondegenerate 
space $V$---its discriminant and its Witt index. 
The discriminant $\delta$ is the determinant of any Gram matrix
for the space. It is either $+1=+$ or $-1=-$ and determines $V$
up to isometry. Concretely, $V$ has discriminant $+1$
if and only if it possesses an orthonormal basis.

The Witt index (introduced in the previous section
for $\gf{2}$-spaces) is the maximum dimension of a totally
singular subspace ($q$ identically $0$). In even dimension $m=2a$
the Witt index is either $a$ or $a-1$, and (again as before)
we attach the Witt sign $\eps$ equal to $+=+1$ or $-=-1$ in 
these respective cases. In odd dimension $m=2a+1$, the Witt index is
always $a$. In even dimension $m$ always
\[
\delta \eps = -1^{\binom{m+1}{2}}\,,
\]
and we use this identity to define the sign $\eps$ for
odd dimension $m$ as well. The space $V$ may then be denoted
$\presub{\delta}V^\eps_m$, which is sometimes abbreviated to
$\presub{\delta}V_m$ or $V^\eps_m$ since, once $m$
has be fixed, the parameters $\delta$ and $\eps$ determine 
each other.

For odd $m$ an equivalent geometric definition
is that $V_m$ has sign $\eps$ when it is isometric to $x^\perp$
for a $+$-vector $x$ in the even dimensional 
$V_{m+1}^{\eps}$.\footnote{Our convention for $\eps$ 
	is that of \cite{CuHa92,CuHa95}. See 
	\cite{CuHa92} for a discussion and a comparison with other
	conventions from
	the literature. Our choice differs from that of 
	Brouwer \cite{aeb} where 
	$\delta \eps = -1^{\binom{m}{2}}$.
	With Brouwer's convention, for odd $m$ and
	$x$ a $+$-vector within $V_m$, the even dimensional
	$x^\perp$ is isometric to $V_{m-1}^{\eps}$.}

As before, in the polar space of $\presub{\delta}V^\eps_m$ there
are exactly two types of $2$-spaces
spanned by singular ($=$ isotropic) vectors:
the totally singular $2$-spaces with $q+1=4$ 
pairwise perpendicular
singular $1$-subspaces and the hyperbolic $2$-spaces
with $s_2=2$ nonperpendicular singular $1$-subspaces
plus a $+$-space and a perpendicular $-$-space. The hyperbolic
$2$-spaces have type $\presub{-}V^+_2$.

The $2$-spaces spanned by nonsingular vectors have three types.
The only nondegenerate example is the asingular space $\presub{+}V^-_2$, 
which is spanned by a pair of perpendicular $+$-spaces and a pair of 
perpendicular $-$-spaces. The two degenerate examples are
the $+$-tangent spaces---consisting of a singular radical of 
dimension $1$ and three $+$-spaces---and the similar $-$-tangent
spaces.

By Witt's Theorem again, the full isometry group of $\presub{\delta}V^\eps_m$
has rank $3$ on the $+$-spaces.\footnote{Similar remarks hold for the
	$-$-spaces, but this only leads to examples isomorphic to the ones
	being discussed.} 
The reflections with centers of $+$-type
form a normal set of $3$-transpositions, commuting pairs of reflections 
corresponding to asingular $2$-spaces $\presub{+}V^-_2$ and noncommuting
pairs to $+$-tangent spaces and their three pairwise nonperpendicular
$+$-spaces.
We consider the $3$-transposition groups $(G,D)$:
the group $G=\omgp[\delta]{m}{\eps}=\omgp{m}{\eps}=\omgp[\delta]{m}{}$ 
(with $m \ge 4$)
is that subgroup of the full isometry group generated by the 
reflection class $D$ having centers of $+$-type. 

\begin{prop}
    For $m \ge 1$,
    the number of singular $1$-spaces
    in $\presub{\delta}V^\eps_m$ is
    \[
    \presub{\delta}s^\eps_m =
        \left\{
        \begin{array}{ll}
        \onehalf(3^{m-1} -1) &\text{for }m\text{ odd};\\
        \onehalf(3^{m-1} -1) + \eps 3^{(m-2)/2} &\text{for }m\text{ even}.
        \end{array}
        \right.
    \]
\end{prop}

\pf
As in \eqref{eqn-s-recur}, the decomposition $V_m^\eps = V^+_2 \perp V_{m-2}^\eps$ implies
\[
s^\eps_m  
= 1 + 3^{m-2} + 3 s^\eps_{m-2} \ .  
\]
When initialized with
\[
s_1^\eps=0\,,\
s_2^+=2\,,\
s_2^-=0\,,
\]
the result follows. \qed

\begin{thm}\label{thm-I5}

\mbox{}

\begin{enumerate}[{\rm(a)}]
\item
For odd $m \ge 5$ the codiagram $\codiag{\omgp{m}{\eps}}$ has extended
parameters
\begin{multline*}
	\bigl(
	(3^{m-1}-\eps3^{(m-1)/2})/2\,,\,
	k'=(3^{m-2}+\eps3^{(m-3)/2})/2,\\
	\lam'=
	\mu'=(3^{m-3}+\eps3^{(m-3)/2})/2
	\,;\\
	\eigmult{r'}{g}=\eigmult{3^{(m-3)/2}}{g}\,,\,
	\eigmult{s'}{f}=\eigmult{-3^{(m-3)/2}}{f}
	\bigr)\,,
\end{multline*}
and
the diagram $\diag{\omgp{m}{\eps}}$ has extended
parameters
\begin{multline*}
	\bigl(
	(3^{m-1}-\eps3^{(m-1)/2})/2\,,\,
	k=3^{m-2}-2\eps3^{(m-3)/2}-1\,,\\
	\lam=2(3^{m-3}-\eps3^{(m-3)/2}-1)\,,\,
	\mu=2(3^{m-3}-\eps3^{(m-3)/2})
	\,;\\
	\eigmult{r}{f}=\eigmult{3^{(m-3)/2}-1}{f}\,,\,
	\eigmult{s}{g}=\eigmult{-3^{(m-3)/2}-1}{g}
	\bigr)\,,
\end{multline*}
where
\[
f=(3^{m-1}-1-(\eps-1)(3^{(m-1)/2}-1))/4
\]
and
\[
g=(3^{m-1}-1-(\eps+1)(3^{(m-1)/2}+1))/4\,.
\]
\item 
For even $m \ge 4$
the codiagram $\codiag{\omgp{m}{\eps}}$ has extended
parameters
\begin{multline*}
	\bigl(
	(3^{m-1}-\eps3^{(m-2)/2})/2\,,\,
	k'=(3^{m-2}-\eps3^{(m-2)/2})/2\,,\\
	\lam'=(3^{m-3}+\eps3^{(m-4)/2})/2\,,\,
	\mu'=(3^{m-3}-\eps3^{(m-2)/2})/2
	\,;\\
	\{
	\eigmult{r'}{g},
	\eigmult{s'}{f}
	\}=
	\{\,
	\eigmult{\eps3^{(m-2)/2}}{d}\,,\,\eigmult{-\eps3^{(m-4)/2}}{e}
	\,\}
	\bigr)\,,
\end{multline*}
and
the diagram $\diag{\omgp{m}{\eps}}$ has extended
parameters
\begin{multline*}
\bigl(
(3^{m-1}-\eps3^{(m-2)/2})/2\,,\,
k=3^{m-2}-1\,,\\
\lam=2(3^{m-3}-1)\,,\,
\mu=2(3^{m-3}+\eps3^{(m-4)/2})
\,;\\
\{
\eigmult{r}{f},\eigmult{s}{g}
\}=
\{\,
\eigmult{-\eps3^{(m-2)/2}-1}{d}\,,\,\eigmult{\eps3^{(m-4)/2}-1}{e}
\,\}
\bigr)\,,
\end{multline*}
where
\[
d=(3^{m/2}-\eps)(3^{(m-2)/2}-\eps)/8
\]
and
\[
e=(3^m-9)/8\,.
\]
\end{enumerate}
\end{thm}

\pf
Some of the calculations work better in terms of $\delta$ while
for others $\eps$ may be preferred. As $\eps$ is the canonical parameter
in even dimension, we state the final results in terms of it,
remembering that always
\[
\delta \eps = -1^{\binom{m+1}{2}}\,.
\]
Some rules-of-thumb for a fixed $\delta$:
\begin{quote}
if $m$ is even then dropping to $m-1$ does not change $\eps$, while
if $m$ is odd then dropping to $m-1$ changes $\eps$ to $-\eps$;
thus any drop by $2$ changes $\eps$ to $-\eps$;
\end{quote}
\begin{enumerate}[$(i)$]
	\item 
    $\presub{\delta}k^\eps_m = 2\,\presub{\delta}s_{m-1}=
    \left\{
    \begin{array}{ll}
    2s_{m-1}^\eps=3^{m-2} -1 &\text{for } m \text{ even};\\
    2\,s_{m-1}^{-\eps}=
    3^{m-2} -1 - 2 \eps 3^{(m-3)/2} &\text{for } m\text{ odd}.
    \end{array}
    \right.$
		
	In the decomposition
	$\presub{\delta}V_m = \presub{+}V_1 \perp \presub{\delta}V_{m-1}$,
	every $+$-tangent on $\langle x \rangle = \presub{+}V_1$
	is spanned by $x$ and the unique singular $1$-space of the
	tangent, which belongs to $x^\perp$. The remaining two
	$+$-spaces of the tangent are adjacent to $\langle x
	\rangle$ in the diagram.
	\item 
	\begin{align*}
		\presub{\delta}(k')_m^\eps &= \presub{\delta}n_{m-1}\\
		&=n_{m-1}^\eps
			= \onehalf(3^{m-2}-\eps 3^{(m-2)/2})  &\text{for}\ m\ \text{even}\,;\\
		&=n_{m-1}^{-\eps}
		= \onehalf(3^{m-2}+\eps 3^{(m-3)/2})
		&\text{for}\ m\ \text{odd}\,;\\
		\presub{\delta}n^\eps_m &=  1 + \presub{\delta}k_m + \presub{\delta}k'_m
		=1+2\,\presub{\delta}s_{m-1}+\presub{\delta}n_{m-1}\\
		&= \onehalf(3^{m-1}-\eps 3^{(m-1)/2})  &\text{for}\ m\ \text{odd}\,;\\
		&= \onehalf(3^{m-1}-\eps 3^{(m-2)/2})
		&\text{for}\ m\ \text{even}\,.
	\end{align*}
	The identity $\presub{\delta}k'_m = \presub{\delta}n_{m-1}$
	follows directly from $\presub{\delta}V_m = \presub{+}V_1 \perp 
	\presub{\delta}V_{m-1}$. Initialization of the recursion is
	provided by 
	\[
	\presub{-}n_1^+=0\,,\
	\presub{+}n_1^-=1\,,\
	\presub{-}n_2^+=1\,,\
	\presub{+}n_2^-=2\,.
	\]
	\item 
	\begin{align*}
		\presub{\delta}(\lam')_m^\eps &= \presub{\delta}n_{m-2}=
		n_{m-2}^{-\eps}\\
		&
		= \onehalf(3^{m-3}+\eps 3^{(m-3)/2})  &\text{for}\ m\ \text{odd}\,;\\
		&
		= \onehalf(3^{m-3}+\eps 3^{(m-4)/2})
		&\text{for}\ m\ \text{even}\,.
	\end{align*}
		The identity $\presub{\delta}(\lam')_m^\eps = \presub{\delta}n_{m-2}$
		follows directly from $\presub{\delta}V_m = \presub{+}V_2 \perp 
		\presub{\delta}V_{m-2}$. 
	\item The parameters we have found so far are
	enough to calculate all remaining ones using the identities of
	Section \ref{sec-rk3}. Some are also geometrically evident.
	Consider the decomposition
	\[
	\presub{\delta}V_m=\presub{+}V_1 \perp \presub{-}V^+_2 \perp
	\presub{-\delta}V_{m-3}\,.
	\]
	Let $x$ be a $+$-vector spanning $\presub{+}V_1$. If $z$
	is a nonzero singular vector in the hyperbolic $\presub{-}V^+_2$ 
	then the $2$-space
	$\langle x,z \rangle$ is a $+$-tangent, and within it $y=x+z$
	spans a $+$-space not perpendicular to $x$. This leads to
	\[
	\presub{\delta}\mu'_m = 3\, \presub{-\delta}n_{m-3}
	\quad\text{and}\quad
	\presub{\delta}\lam_m = \presub{\delta}s_{m-1} 
	+ 3\, \presub{-\delta}s_{m-3}\,.
	\]
	\qed
\end{enumerate}

\begin{prop}\label{prop-P5-case}

\mbox{}

	\begin{enumerate}[$(a)$]
		\item
			\name{PR5:} the diagram 
			$\diag{3^\comd{h}\,\omgp{m}{\eps}}$
			for odd $m \ge 5$,
			 $\eps=\pm$,
			and $h \ge 0$
			has size 
        	\[
            n=3^h(3^{m-1}-\eps3^{(m-1)/2})/2
            \] 
            and spectrum
            \[
            \spec{3^{m-2+h}-2\eps3^{(m-3)/2+h}-1\,;\,
            \eigmult{3^{(m-3)/2+h}-1}{f}\,,\,
            [-1]^\star,\,
            \eigmult{-3^{(m-3)/2+h}-1}{g}
            }
            \]
			where
			\[
			f=(3^{m-1}-1-(\eps-1)(3^{(m-1)/2}-1))/4
			\]
			and
			\[
			g=(3^{m-1}-1-(\eps+1)(3^{(m-1)/2}+1))/4\,.
			\]
		\item
			\name{PR5:} the diagram 
			$\diag{3^\comd{h}\,\omgp{m}{\eps}}$
			for even $m \ge 6$,
			$\eps=\pm$,
			and $h \ge 0$
			has size  
            \[
            n=3^h(3^{m-1}-\eps3^{(m-2)/2})/2
            \] 
            and spectrum
            \[
            \spec{3^{m-2+h}-1\,
            ;\,\eigmult{-\eps3^{(m-2)/2+h}-1}{d}\,,\,
            [-1]^\star,\,
            \eigmult{\eps3^{(m-4)/2+h}-1}{e}}
            \]
			where
			\[
			d=(3^{m/2}-\eps)(3^{(m-2)/2}-\eps)/8
			\]
			and
			\[
			e=(3^m-9)/8\,.
			\]
		\end{enumerate}
	\qed
\end{prop}

\begin{prop}\label{prop-P8-case}
	\name{PR8:} the diagram 
	$\diag{4^\comd{h}\splt(\osix)}$ for $h \ge 0$
	has size 
	\[
	n=126(4^{h})
	\] 
	and spectrum
	\[
	\spec{80(4^{h})\,;\eigmult{8(4^{h})}{35},
		\eigmult{-4^{h+1}}{90},\eigmult{0}{\star}}\,.
	\]
\end{prop}

\pf
Apply Corollary \ref{cor-step} to the
diagram $\diag{\omgp{6}{-}} = \diag{\osix}$,
which has extended parameters
\[
\left(126,80,52,48\,; 
\eigmult{8}{35},\eigmult{-4}{90}
\right)\,.
\]
\qed

\subsection{Sporadic cases}

\begin{thm}\label{thm-PR7}
\begin{enumerate}[{\rm(a)}]
	\item The codiagram $\codiag{\Fi{22}}$ has extended
	parameters
	\[
	\left(3510,693,180,126\,; 
	\eigmult{63}{429},
	\eigmult{-9}{3080}
	\right)
	\,.
	\]
  \item \name{PR7(a):} 
  The diagram $\diag{\Fi{22}}$ has extended
parameters
\[
  \left(3510,2816,2248,2304\,; 
  \eigmult{8}{3080},\eigmult{-64}{429}
\right)
\]
and spectrum
\[
\spec{2816\,; 
\eigmult{8}{3080},\eigmult{-64}{429}}\,.
\]
  \item The codiagram $\codiag{\Fi{23}}$ has extended
parameters
\[
  \left(31671,3510,693,351\,; 
  \eigmult{351}{782},\eigmult{-9}{30888}
  \right)
\,.
\]
  \item \name{PR7(b):} 
  The diagram $\diag{\Fi{23}}$ has extended
  parameters
  \[
  \left(31671,28160,25000,25344\,; 
  \eigmult{8}{30888},\eigmult{-352}{782}\right)
  \]
  and spectrum
  \[
  \spec{
  	28160\,;
  	\eigmult{8}{30888},\eigmult{-352}{782}
  }
  \]
  \item The codiagram $\codiag{\Fi{24}}$ has extended
parameters
\[
  \left(306936,31671,3510,3240\,; 
  \eigmult{351}{57477}
  \eigmult{-81}{249458},
  \right)
\,.
\]
  \item \name{PR7(c):} 
  The diagram $\diag{\Fi{24}}$ has extended
  parameters
  \[
  \left(306936,275264,246832,247104\,; 
  \eigmult{80}{249458},\eigmult{-352}{57477}\right)
  \]
  and spectrum
  \[
  \spec{275264\,; 
  \eigmult{80}{249458},\eigmult{-352}{57477}}
  \,.
  \]
\end{enumerate}
\end{thm}

\pf 
See \cite{Fi71} or \cite{Hu75} for the basic
parameters.  The extended parameters can then be
calculated as in Section \ref{sec-rk3} and
are also given in \cite{aeb}.

While we do not repeat these calculations,
the basic parameters for the codiagram ($=$ commuting graph)
appear naturally within Fischer's $3$-transposition theory.
Fischer \cite{Fi71} attacked the classification
by induction, noting
that the $3$-transposition group $(G,D)$
is essentially determined by the codiagrams of two 
if its ``local'' $3$-transposition subgroups:
\[
K_G= \langle \oC_D(d) \rangle
\quad\text{and}\quad
M_G= \langle \oC_D(d,c) \rangle
\]
for $d,c \in D$ with $|dc|=3$. 
Fischer used this local data to reconstruct the
global group $(G,D)$.

This is particularly relevant for us since
\[
k'= |\oC_D(d) \backslash \{d\}|=|\codiag{K_G}|
\]
and
\[
\mu'= |\oC_D(d,c)|=|\codiag{M_G}|\,.
\]
The additional $3$-transposition subgroup 
$L_G=\langle \oC_D(d,e) \rangle$  
with $|de|=2$ 
(naturally found as the subgroup $K_{K_G}$
of $K_G$)
yields
\[
\lam'= |\oC_D(d,e) \backslash \{d,e\}|=|\codiag{L_G}|\,.
\]
The local parameters $k',\mu',\lam'$ then allow us,
using \eqref{eqn-mul}, to
calculate the global parameter
\[
n=1+k'+l'=1+k'+k'(k'-1-\lam')/\mu'\,.
\]

For each pair 
of $3$-transposition group $\tilde{K}$ supplied
with $3$-transposition subgroup $\tilde{M}$,
Fischer looked for a $3$-transposition group $G$
with $\codiag{K_G} = \codiag{\tilde{K}}$ and 
$\codiag{M_G} = \codiag{\tilde{M}}$.
Everything went smoothly, producing the symmetric
and classical groups which we have discussed;
but there was one loose-end---the 
pair 
\[
{\oP\!\su{6}} = \tilde{K} \ge \tilde{M}= \omgp{6}{-}\,.
\]
This special case
led to a tower of $3$-transposition 
groups---those that are sporadic. Specifically
Fischer found the $(\codiag{K},\codiag{M})$-tower:
\begin{multline*}
(\codiag{\oP\!\su{6}},\codiag{\omgp{6}{-}})
=(\codiag{K_{\Fi{22}}},\codiag{M_{\Fi{22}}})\,,\ \\
(\codiag{\Fi{22}},\codiag{\omgp{7}{+}})=(\codiag{K_{\Fi{23}}},
\codiag{M_{\Fi{23}}})\,,\
\\
(\codiag{\Fi{23}},\codiag{\tri{3}})=(\codiag{K_{\Fi{24}}},
\codiag{M_{\Fi{24}}})\,.
\end{multline*}

This construction of the
sporadic examples aids in the identification
and calculation of their basic parameters. In the table,
$\codiag{K_G}$ occurs on the line above $\codiag{G}$
and $\codiag{L_G}$ two lines
above. We have:
\[
\begin{array}{crrrrrrr|rc}
 \codiag{G}  && \{n\,, & k'\,, & \lam'\}&& &&\mu'&\codiag{M_G} \\
 \hline
\codiag{\oP\!\su{6}}	&&&&&\mathit{693}&180& 
51
&45&\codiag{\su{4}}\\
\codiag{\Fi{22}}&	&&&\mathit{3510}&693&180&&126&\codiag{\omgp{6}{-}} \\
\codiag{\Fi{23}}&	&&\mathit{31671}&3510&693&&&351&\codiag{\omgp{7}{+}}\\
\codiag{\Fi{24}}&	&\mathit{306936}&31671&3510&&&&3240&\codiag{\tri{3}}\\
\end{array}
\]
Here the global parameters (in italics) can be calculated from the local
parameters. Initialization is provided by the values for
$\codiag{\su{6}}=\codiag{\oP\!\su{6}}$ found
in Theorem \ref{thm-I6}\footnote{Although the codiagrams
	$\codiag{\oP\!\su{6}}$ and $\codiag{\su{6}}$
	are equal, we use the first notation here because
	$K_{\Fi{22}}/\langle d \rangle$
	is isomorphic to $\oP\!\su{6}$.
	Recall that central type is a coarser
	equivalence relation than isomorphism.
}
; the
remaining $M_G$ were identified as part of Fischer's induction. 
For the size $n$ of $\codiag{\Fi{22}}$ we calculate
\[
1+693+693(693-1-180)/126 = 3510
\]
as claimed. The others are similar.\qed

\begin{thm}\label{thm-tri}
	\begin{enumerate}[{\rm(a)}]
		\item 
			\name{PR7(d):} the diagram 
			$\diag{\tri{2}}$,
			has size 
			\[
			n=360
			\] 
			and spectrum
			\[
			\spec{296;[-64]^2,[8]^{105},[-4]^{252}}\,.
			\]
		\item
		 	\name{PR7(e):} 
		 	the diagram 
		 	$\diag{\tri{3}}$,
		 	has size 
		 	\[
		 	n=3240
		 	\] 
		 	and spectrum
		 	\[
		 	\spec{2888;[-352]^2,[8]^{2457},[-28]^{780}}\,.
		 	\]
	\end{enumerate}
\end{thm}

\pf
By Theorem \ref{thm-I3}
the diagram $\diag{\ogpt[+]{8}}=\diag{\oP\!\Omega_8^+(2)\splt\Sym{2}}$ 
has extended parameters
 \[
 \left(120,56,28,24\,; 
 \eigmult{8}{35},\eigmult{-4}{84}
 \right)\,.
 \]
As in Proposition \ref{prop-eigtri}, (the adjacency matrix for) the
diagram $\diag{\tri{2}}$ is $3 \times \diag{\ogpt[+]{8}}$,
so it has size $3(120)=360$ and spectrum 
\begin{multline*}
\spec{56+2(120);\eigmult{-(120-56)}{2},
	\eigmult{8}{3(35)},
	\eigmult{-4}{3(84)}}=\\
\spec{296;[-64]^2,[8]^{105},[-4]^{252}}\,.
\end{multline*}

Similarly 
$\diag{\tri{3}}$ is $3 \times \diag{\omgp{8}{+}}$,
where the diagram $\diag{\omgp{8}{+}}$ 
$(=\diag{\oP\!\Omega_8^+(3)\splt\Sym{2}})$ has, 
by Theorem \ref{thm-I5}, extended parameters 
\[
 \left(1080,728,484,504\,; 
 \eigmult{8}{819},\eigmult{-28}{260}
 \right)\,.
\]

Proposition \ref{prop-eigtri} again applies to give the
spectrum
\begin{multline*}
	\spec{728+2(1080);\eigmult{-(1080-728)}{2},
			\eigmult{8}{3(819)},
			\eigmult{-28}{3(260)}
		}
	=\\
	\spec{2888;[-352]^2,[8]^{2457},
		[-28]^{780}}\,.
\end{multline*}
\qed

\section{Diagram minimum eigenvalues}

Miyamoto \cite{Mi96} associated
$3$-transposition groups with the Griess algebras
of certain vertex operator algebras of $OZ$-type.
In that context, the minimum eigenvalue of the
diagram for the group is important, particularly
those with minimum eigenvalue greater than
or equal to $-8$. Classification
of the associated groups and Griess algebras
was pursued by Miyamoto and Kitazume \cite{KM01}
and Matsuo \cite{Ma03,Ma05}. Similar issues
arise for minimum eigenvalue at least $-64$, and that
was the initial motivation for the current
paper.

\subsection{Compact Matsuo and Griess algebras}

Let $\eta$ be an element of $\bbR$ not equal to $0$ or $1$.
A real \announce{Matsuo algebra} for the eigenvalue $\eta$ 
is a commutative algebra $M=\bigoplus_{a \in \mathcal{A}} \bbR a$
with basis $\mathcal{A}=\{ a_i \mid i \in I\}$ of idempotents
$a_i$ (called \announce{axes}) and having the property that any two
$a,b \in \mathcal{A}$ generate
one of the subalgebras
\begin{enumerate}[$(i)$]
	\item $1A=\bbR$ with $a=b$;
	\item $2B=\bbR^2= \bbR a \oplus \bbR b$ with $ab=0$;
	\item $3C(\eta) = \bbR a \oplus \bbR b \oplus \bbR c$
	with $xy=\frac{\eta}{2}(x+y-z)$ for $\{x,y,z\}
	=\{a,b,c\} = (\bbR a \oplus \bbR b \oplus \bbR c) \cap \mathcal{A}$.
\end{enumerate}
On $M$ we define the symmetric bilinear form 
$\lla \cdot \!\mid\! \cdot \rra$
given by, respectively,
\begin{align*}
	1A:\quad \lla a\!\mid\! a \rra &= 1\,;  \\
    2B:\quad \lla a\!\mid\! b \rra  &=0\,; \\
    3C(\eta):\quad \lla a\!\mid\! b \rra &= \frac{\eta}{2}\,.
\end{align*}
We say that the algebra $M$ is \announce{compact} if the
associated form is positive definite. 

Matsuo algebras were introduced
\cite{MM99,Ma03,Ma05} because certain compact Matsuo
algebras arise as the Griess algebras of compact
vertex operator algebras of $OZ$ type, as noted
by Miyamoto \cite{Mi96}. A classification of all such Griess
algebras is desirable.

The crucial observation, due to Miyamoto, is
that in the Griess algebra case, for each axis
$a \in \mathcal{A}$, the permutation  $\tau_a$
of $\mathcal{A}$ given by
\begin{align*}
	1A:\quad  &\tau_a(a) = a\,;  \\
	2B:\quad  &\tau_a(b) = b\,;  \\
	3C(\eta):\quad &\tau_a(b) = c\,.
\end{align*}
is an automorphism and indeed
$\{\tau_{a_i} \mid i \in I\}$ is
a normal set of $3$-transpositions in the automorphism
group of $M$. It is enough to consider the case
in which $D=\{\tau_{a_i} \mid i \in I\}$ is
a class of $3$-transpositions in the group
$G=\langle D \rangle \le \operatorname{Aut}(M)$.

This property of Griess and Matsuo algebras was
seen in \cite{HRS15} to characterize\footnote{This
	remark is somewhat inaccurate in the special
	case $\eta=\frac{1}{2}$. See \cite{HRS15,HSS18}
	for precision.}
the \announce{axial algebras of Jordan type $\eta$}.
In that case the symmetric form
$\lla \cdot \!\mid\! \cdot \rra$
is \announce{associative} in that
\[
\lla xy \!\mid\! z \rra = \lla x \!\mid\! yz \rra
\]
for all $x,y,z \in M$. An important property of
every associative form is that its radical $R$ is
an ideal of the algebra $M$.

The \announce{Gram matrix} of the form with
respect to $\mathcal A$ is
\[
I + \frac{\eta}{2}H\,,
\]
where $H$ is the adjacency matrix of the diagram $\diag{D}$,
so the compact axial algebra $M$ of Jordan type $\eta$ 
must have 
\[
1 + \frac{\eta}{2}\rho > 0
\]
hence
\[
\rho > -\frac{2}{\eta}\,,
\]
where $\rho$ is the minimum eigenvalue of $\diag{D}$.
In the case $\rho = -\frac{2}{\eta}$, the algebra
$M$ is positive semidefinite and its axial
quotient $M/R$ is again compact (which is to say, 
positive definite). The quotient is thus also a 
candidate to be a compact Griess algebra.

Initial interest focuses on the eigenvalues 
\[
\eta = \frac{1}{4} \quad\text{and}\quad \eta=\frac{1}{32}
\]
since these are the eigenvalues associated with the
Monster algebra of Griess \cite{Gr81} as embedded
in the Moonshine vertex operator algebra of
Frenkel, Lepowsky, and Meurman \cite{FLM88}.
Correspondingly, we are interested in the minimum
eigenvalues 
\[
\rho \ge -8 \quad\text{and}\quad \rho \ge -64\,.
\]

The Griess algebra case for the eigenvalue $\eta=\frac{1}{4}$ 
was investigated by Kitazume and Miyamoto \cite{KM01}
and Matsuo \cite{Ma03,Ma05}. The classification
is due to Matsuo:

\begin{thm}\label{thm-matsuo-eig}
	A compact Griess algebra for the eigenvalue $\eta=\frac{1}{4}$
	exists if and only if the associated finite $3$-transposition
	group has one of the central types below.
	In each case the algebra is uniquely determined
	as the corresponding Matsuo algebra 
	(modulo the radical of its form when the minimum
	eigenvalue is $\rho=-8$).
	\begin{enumerate}[$(a)$]
		\item 
		$\rho=-1:$ \quad
		$\Sym{3}$;
		\item 
		$\rho=-2^{r+1}$ with $0 \le r \le 2:$
		\quad
		$(2^{m-1})^r\splt\Sym{m}$
		for each $m\ge 4$;
		\item 
		one of the nine individual groups
		\begin{enumerate}[$(i)$]
			\item $\rho=-4:$\quad
			$\ogpt[-]{6};$\
			$\ogpt[+]{8};$\
			$\spt{6};$\
			\item $\rho=-8:$\quad
			$2^6\splt\ogpt[-]{6};$\
			$\ogpt[-]{8};$\
			$2^8\splt\ogpt[+]{8}$;\
			$\ogpt[+]{10};$\
			$2^6\splt\spt{6};$\
			$\spt{8}$.
		\end{enumerate}
	\end{enumerate}
\end{thm}

Matsuo's proof comes in three pieces:
\begin{enumerate}[$(i)$]
\item Properties of Griess algebras prove that $(G,D)$
must have the central type of a $3$-transposition
subgroup of some $\spt{2n}$. (See \cite[Prop.~1]{Ma05}.)
\item Identification of the $3$-transposition
groups of symplectic type with minimum eigenvalue
$\rho \ge -8$. (These are precisely the groups of
the theorem, as can be verified from
the table in Section \ref{sec-eig}
or the results of the Subsection \ref{sec-64}.)
\item For each qualifying group, checking that the appropriate
Matsuo algebra (quotient) is indeed a Griess algebra
for some vertex operator algebra.
\end{enumerate}

\subsection{Classification by minimum eigenvalue}

The second table of Section \ref{sec-eig} provides the
minimum eigenvalue in bold for
the diagram of each finite $3$-transposition groups
and so has several useful consequences. In parallel to
Matsuo's result Theorem \ref{thm-matsuo-eig} we have:

\begin{thm}\label{thm-minl-eig}
	There are nondecreasing, nonnegative integral valued 
	functions $S(t)$ and $I(t)$ defined on $2 \le t \in \bbZ^+$ 
	such that the diagram $\diag{D}$ of a finite conjugacy class 
	$D$ of $3$-transpositions has minimum eigenvalue $\rho_{\min} \ge -t$ 
	if and only if the corresponding $3$-transposition group $(G,D)$ 
	belongs to one of the following central type classes: 
	\begin{enumerate}[$(a)$]
		\item 
		infinitely many groups $3^u \splt 2$ of Moufang type
		(the case $\rho_{\min}=-1$);
		\item 
		$S(t)$ distinct	groups
		{$N \splt \Sym{m}$} for each $m\ge 4$;
		\item 
		$I(t)$  {individual examples}. \qed
	\end{enumerate}
\end{thm}

For a given $\eta \le \frac{2}{t}$
let $S^\eta(t)$ and $I^\eta(t)$ be the corresponding functions 
counting those $3$-transposition groups realized by
some Griess algebra for the eigenvalue $\eta$. Clearly
$0 \le S^\eta(t) \le S(t)$ and $0 \le I^\eta(t) \le I(t)$.
Matsuo's Theorem \ref{thm-matsuo-eig} and the results
of the next subsection give
\[
3=S^{\frac{1}{4}}(8) \le S(8)=4
\quad\text{and}\quad
9=I^{\frac{1}{4}}(8) \le I(8)=14\,.
\]
The differences are caused by Matsuo's proven restriction
to symplectic type. For $\rho \ge -8$ 
three of the four
symmetric families as in Theorem \ref{thm-minl-eig}(b)
have symplectic type ($\Weyl[3]{\tilde{A}_{m-1}}$ does not),
and Matsuo showed that all lead to Griess algebras.
Similarly exactly nine of the $14$ individual groups
counted by $I(8)$ are of symplectic type, and they
too produce Griess algebras. Almost by definition,
$\Sym{3}$ is the only $3$-transposition group that
is simultaneously of Moufang and symplectic type,
and this is reflected in the stark difference between
Theorem \ref{thm-matsuo-eig}(a) and
Theorem \ref{thm-minl-eig}(a).

For the case $\eta=\frac{1}{32}$, hence $\rho \ge -64$,
the next section reveals
\[
S(64)=13
\quad\text{and}\quad
I(64)=90\,.
\]
Very little is known about the corresponding
$S^{\frac{1}{32}}(64)$ and $I^{\frac{1}{32}}(64)$.
Chen and Lam \cite{ChLa14} have shown that $\su{3}'$ can be
realized for $\eta=\frac{1}{32}$. In particular
Matsuo's restriction to symplectic type will
not be available in this case.

\subsection{Minimum eigenvalue $\rho \ge -64$}
\label{sec-64}

The table also yields:
\begin{thm}
	\label{thm-four-cases}
	Let $(G,D)$ be a finite $3$-transposition group.
	Then the minimum eigenvalue $\rho$ of its diagram $\diag{G}$
	satisfies one of:
	\begin{enumerate}[$(a)$]
		\item
		$\rho=-1$ and $G$ has Moufang type $3^u \splt 2$.
		\item 
		$\rho=-2^{a}$ for the positive integer $a$, 
		with $(G,D)$ being one of:
		\begin{enumerate}[$(i)$]
		\item an infinite class of examples with quotient $\Sym{m}$ 
		under \name{PR2(a)};
		\item if $a$ is even and at least $4$, an infinite class of examples 
		with quotient $\Sym{m}$ under \name{PR2(d)};
		\item a finite number of classical examples 
		in characteristic $2$ under \name{PR3,4,6,17-19};
		\item if $a$ is even, a single mixed characteristic example
		$(4^6)^h\splt\osix$ with $h=(a-2)/2$ under \name{PR5,8};
		\item if $a=6$ (so $\rho=-64$), the examples $\Fi{22}$ or $\tri{2}$ 
		under \name{PR7}.
	\end{enumerate}			
		\item 
		$\rho=-3^{b}-1$ for the positive integer $b$, 
		with $(G,D)$ being one of:
		\begin{enumerate}[$(i)$]
			\item an infinite class of examples with quotient $\Sym{m}$ 
			under \name{PR2(b)};
			\item if $b$ is at least $2$, an infinite class of examples 
			with quotient $\Sym{m}$ under \name{PR2(c)};
			\item a finite number of classical examples 
			in characteristic $3$ under \name{PR5,13-16};
			\item a finite number of mixed characteristic examples
			under \name{PR9-12}.
		\end{enumerate}
		\item
		$\rho=-352$ and $(G,D)$ has type $\Fi{23}$, $\Fi{24}$,
		or $\tri{3}$ under \name{PR7}. \qed
	\end{enumerate}
\end{thm}

Here we list those groups whose diagrams have minimum eigenvalue
greater than or equal to $-64$. 
As $3^b+1$ is never a multiple of $8$, cases (a) and
	(b) only overlap at $\rho=-4=-2^2=-3^1-1$. This case
	is discussed in Section \ref{subsec-r2}.

We list the actual groups---as detailed in Section 2 of 
\cite{CuHa95}---not
just their diagrams
(although we continue not to distinguish between groups
of the same central type unless necessary). 
Therefore, unlike Theorem
\ref{thm-cuha}, the exotic cases \name{PR13-19} 
from \cite{CuHa95} 
appear separately from other cases. 
These all involve nonsplitting
of certain extensions, and their diagrams are the same
as those of groups from earlier in the list all of whose extensions
split. For instance, the split extension $4^7 \splt \su{7}$
of \name{PR6}\ and the nonsplit extension $4^7 \nsplt \su{7}$
of \name{PR18}\ have same diagram, and so they share the 
minimum eigenvalue $-64$.
Similarly for $\rho=-28$ the groups
$(3^5)^2\splt\omgp{5}{-}$ of \name{PR5} 
and
$(3^5\odot 3^5)\splt\omgp{5}{-}$ 
of \name{PR13} 
have the same diagram but differ in that each
is the split extension of a subgroup $3^{10}$ 
by $\omgp{5}{-}$, but the first $3^{10}$
is, as $\omgp{5}{-}$-module,
the direct sum 
$(3^5)^2=3^5\oplus 3^5$
of two
copies of the natural module $3^5$ while the second is 
a nonsplit $\omgp{5}{-}$-module extension 
$3^5\odot 3^5$.
(Here we use notation where
a split module extension with submodule $A$ and quotient $B$ 
is denoted $A \oplus B$ while a nonsplit module extension is 
$A\odot B$.

The four classes with symmetric quotient mentioned in parts (b) and (c) of the theorem 
are those of PR2(a), PR2(b), PR2(c), and PR2(d). For the classes  
PR2(a), PR2(c), and PR2(d) the parameters $a$ or $b$ and $m$ 
determine the group $G$ uniquely up to central type. That is false
for PR2(b) with central type $3^\comd{h}\splt\Sym{m}$ whenever
$h \ge 3$. In that case, the type of $G$ is that of some
$\operatorname{Wr}(B,m)$---the subgroup of the wreathed product 
of $B$ by $\Sym{m}$ generated by its transpositions
\cite[p.\,162]{CuHa95}. Here $B$ can be any group of
exponent $3$. For $|B|=3^b$ the group $\operatorname{Wr}(B,m)$ 
is a $3$-transpostion group of central type
$3^\comd{b}\splt\Sym{m}$ having minimum eigenvalue $-3^b-1$. As $b$ increases,
the number of choices for $B$ increases dramatically.
The smallest nontrivial case is $b=3$ where the only two choices for 
$B$ are the elementary abelian group $3^3$ and the extraspecial group
$3^{1+2}$. This leads to two different groups with minimum eigenvalue
$-28$, as seen in Section \ref{ss-28} below.

\subsubsection{$\rho=-1$}

\mbox{}

\name{PR1.} 
The examples are the groups $3^u \splt 2$ of Moufang type
with complete diagram. (Indeed, any connected regular graph with
minimum eigenvalue $-1$ is complete. Exercise!)

\subsubsection{$\rho=-2$}
\label{subsec-r2}
\mbox{}

\name{PR2(a).} 
A $3$-transposition group has minimum eigenvalue $-2$ if and
only if it is isomorphic to $\Sym{m}\,(=\Weyl{A_{m-1}})$ for some $m \ge 4$.

\subsubsection{$\rho=-4$}

The eigenvalue $-4$ is anomalous, as it can
be written $-2^2$, as in Theorem \ref{thm-four-cases}(b),
and $-3^1-1$, as in Theorem \ref{thm-four-cases}(c).
Thus it behaves like a characteristic $2$ case and also
like a characteristic $3$ case. Both parts of the theorem
predict two infinite families with symmetric quotient.
In the characteristic $2$
case these should be \name{PR2(a)}\ and \name{PR2(d)}\
while in the characteristic $3$ case these should be
\name{PR2(b)}\ and \name{PR2(c)}. The eigenvalue $-4$
compromises by choosing \name{PR2(a)}\ and \name{PR2(b)}.
We also have the mixed characteristic example $\omgp{6}{-}$.

Another mixed characteristic oddity for $-4$ is that
the groups
\[
\omgp{5}{+} =\ogpt[-]{6}
\quad \text{and}\quad
\omgp{5}{-} = 2 \times \su{4}
\]
appear twice on the list, once under \name{PR5}\
in characteristic $3$
and a second time under \name{PR3}\ or \name{PR6},
as appropriate, in characteristic $2$.

\dspace{1}

\name{PR2(a).} 
$2^{m-1}\splt \Sym{m}\,(=\Weyl{D_m}=\Weyl[2]{\tilde{A}_{m-1}})$
for all $m \ge 4$.

\name{PR2(b).} 
$3^{m-1}\splt\Sym{m}\,(=\Weyl[3]{\tilde{A}_{m-1}})$ for all $m\ge 4$.

\name{PR3.} 
$\ogpt[-]{6}\,(= \Weyl{E_6} = \omgp{5}{+})$;
$\ogpt[+]{8}\,(=\Weyl{E_8}/2)$. 

\name{PR4.}  
$\spt{6}\,(=\Weyl{E_7}/2)$. 

\name{PR5.} 
$\omgp{5}{+}\,(= \Weyl{E_6}=\ogpt[-]{6})$;
$\omgp{5}{-}\,(= 2 \times \su{4})$;
$\omgp{6}{-}$.

\name{PR6.} 
$4^3\splt\su{3}'$;
$\su{4}\,(=\omgp{5}{-}/2)$;
$\su{5}$.

\subsubsection{$\rho=-8$}

\mbox{}

\name{PR2(a).} 
$(2^{m-1})^2\splt\Sym{m}$
for all $m\ge 4$.

\name{PR3.} 
$2^6\splt\ogpt[-]{6}\,(=\Weyl[2]{\tilde{E}_6})$;
$2^8\splt\ogpt[+]{8}\,(=\Weyl[2]{\tilde{E}_8}/2)$;
$\ogpt[-]{8}$;
$\ogpt[+]{10}$.

\name{PR4.} 
$2^6\splt\spt{6}\,(=\Weyl[2]{\tilde{E}_7}/2)$;
$\spt{8}$.

\subsubsection{$\rho=-10$}

\mbox{}

\name{PR2(b).} 
$(3^{m-1})^2\splt\Sym{m}$
for all $m\ge 4$.

\name{PR2(c).} 
$3^m\splt2^{m-1}\splt\Sym{m}\,(=\Weyl[3]{\tilde{D}_m})$
for all $m\ge 4$.

\name{PR5.} 
$3^5\splt\omgp{5}{-}
\,(=\Weyl[3]{\tilde{E}_6}/3)$;\,
$3^5\splt\omgp{5}{+}$;\,
$3^6\splt\omgp{6}{-}$;\,
$\omgp{6}{+}$;\,
$\omgp{7}{-}$;\,
$\omgp{7}{+}$;\,
$\omgp{8}{-}$.

\name{PR9.} 
$3^7\splt(2 \times \spt{6})\,(=\Weyl[3]{\tilde{E}_7}).$

\name{PR10.} 
$3^8\splt(2\nsplt\ogpt[+]{8})\,(=\Weyl[3]{\tilde{E}_8}).$

\subsubsection{$\rho=-16$}

\mbox{}

\name{PR2(a).} 
$(2^{m-1})^3\splt\Sym{m}$
for all $m\ge 4$.

\name{PR2(d).} 
$4^m\splt3^{m-1}\splt\Sym{m}$, 
for all $m\ge 4$.

\name{PR3.} 
$(2^6)^2\splt\ogpt[-]{6}$;\,
$(2^8)^2\splt\ogpt[+]{8}$;\,
$2^8\splt\ogpt[-]{8}$;\,
$2^{10}\splt\ogpt[+]{10}$;\,
$\ogpt[-]{10}$;\,
$\ogpt[+]{12}$.

\name{PR4.} 
$(2^6)^2\splt\spt{6}$;\,
$2^8\splt\spt{8}$;\,
$\spt{10}$.

\name{PR6.} 
$(4^3)^2\splt\su{3}'$;\,
$4^4\splt\su{4}$;\,
$4^5\splt\su{5}$;\,
$\su{6}$;\,
$\su{7}$.

\name{PR8.} 
$4^6\splt(\osix)$.

\subsubsection{$\rho=-28$}
\label{ss-28}

\mbox{}

\name{PR2(b).} 
$(3^{m-1})^3\splt\Sym{m}$
and 
$(3^{m-1}\nsplt(3^{m-1})^2)\splt\Sym{m}$
for all $m\ge 4$.

\name{PR2(c).} 
$(3^{m})^2\splt2^{m-1}\splt\Sym{m}$
for all $m\ge 4$.

\name{PR5.} 
$(3^5)^2\splt\omgp{5}{-}$;\,
$(3^5)^2\splt\omgp{5}{+}$;\,
$(3^6)^2\splt\omgp{6}{-}$;\,
$3^6\splt\omgp{6}{+}$;\,
$3^7\splt\omgp{7}{-}$;\,
$3^7\splt\omgp{7}{+}$;\,
$3^8\splt\omgp{8}{-}$;\,
$\omgp{8}{+}$;\,
$\omgp{9}{-}$;\,
$\omgp{9}{+}$;\,
$\omgp{10}{-}$.

\name{PR9.} 
$(3^6)^2\splt(2 \times \spt{6})$.

\name{PR10.} 
$(3^8)^2\splt(2\nsplt\ogpt[+]{8})$.

\name{PR11.} 
$3^{10}\splt(2 \times \su{5})$.

\name{PR12.}
$3^8\splt (U \splt\su{3}')$, 
$U= 2^{1+6}$, $U'=2$, $U/U'=4^3$.

\name{PR13.} 
$(3^5\odot 3^5)\splt\omgp{5}{-}$.

\name{PR14.} 
$(3^6\odot 3^6)\splt(3\cdot\omgp{6}{-})$.

\name{PR15.} 
$3^7\,\nsplt\omgp{7}{-}$.

\name{PR16.} 
$3^8\nsplt\omgp{8}{-}$.

\subsubsection{$\rho=-32$}

\mbox{}

\name{PR2(a).} 
$(2^{m-1})^4\splt\Sym{m}$
for all $m\ge 4$.

\name{PR3.} 
$(2^6)^3\splt\ogpt[-]{6}$;\,
$(2^8)^3\splt\ogpt[+]{8}$;\,
$(2^8)^2\splt\ogpt[-]{8}$;\,
$(2^{10})^2\splt\ogpt[+]{10}$;\,
$2^{10}\splt\ogpt[-]{10}$;\,
$2^{12}\splt\ogpt[+]{12}$;\,
$\ogpt[-]{12}$;\,
$\ogpt[+]{14}$.

\name{PR4.}  
$(2^6)^3\splt\spt{6}$;\,
$(2^8)^2\splt\spt{8}$;\,
$2^{10}\splt\spt{10}$;\,
$\spt{12}$.

\subsubsection{$\rho=-64$}

\mbox{}

\name{PR2(a).} 
$(2^{m-1})^5\splt\Sym{m}$
for all $m\ge 4$.

\name{PR2(d).} 
$(4^{m})^2\splt3^{m-1}\splt\Sym{m}$, 
for all $m\ge 4$.

\name{PR3.} 
$(2^6)^4\splt\ogpt[-]{6}$;\,
$(2^8)^4\splt\ogpt[+]{8}$;\,
$(2^8)^3\splt\ogpt[-]{8}$;\,
$(2^{10})^3\splt\ogpt[+]{10}$;\,
$(2^{10})^2\splt\ogpt[-]{10}$;\,
$(2^{12})^2\splt\ogpt[+]{12}$;\,
$2^{12}\splt\ogpt[-]{12}$;\,
$2^{14}\splt\ogpt[+]{14}$;\,
$\ogpt[-]{14}$;\,
$\ogpt[+]{16}$.

\name{PR4.}  
$(2^6)^4\splt\spt{6}$;\,
$(2^8)^3\splt\spt{8}$;\,
$(2^{10})^2\splt\spt{10}$;\,
$2^{12}\splt\spt{12}$;\,
$\spt{14}$.

\name{PR6.} 
$(4^3)^3\splt\su{3}'$;\,
$(4^4)^2\splt\su{4}$;\,
$(4^5)^2\splt\su{5}$;\,
$4^6\splt\su{6}$;\,
$4^7\splt\su{7}$;\,
$\su{8}$;\,
$\su{9}$.

\name{PR7.} 
$\Fi{22}; \tri{2}$.

\name{PR8.} 
$(4^6)^2\splt(\osix)$.

\name{PR17.} 
$(4^5 \odot 4^5)\splt\su{5}$.

\name{PR18.} 
$4^7\nsplt\su{7}$.

\name{PR19.} 
$T\splt\su{3}'$,
$T=4^{3+(3+3)}$, 
$\oZ(T)=T'=4^3$,
$T/T'= (4^3)^2$.

\end{document}